\documentclass{article}

\usepackage{arxiv}
\usepackage{cite}
\usepackage{amsmath,amssymb,amsfonts,bm}
\usepackage{mathtools}
\usepackage{algorithmic}
\usepackage{graphicx}
\usepackage[dvipsnames]{xcolor}
\usepackage{textcomp}
\usepackage{dsfont}
\usepackage{url}
\usepackage[outdir=./]{epstopdf}        
\usepackage{enumerate}
\usepackage{lscape}
\usepackage{multicol}
\usepackage{multirow}
\usepackage{calligra}
\usepackage{booktabs}
\usepackage{mathtools}
\usepackage{framed} 
\usepackage{empheq}
\usepackage{graphics} 
\usepackage{epsfig} 
\usepackage{times} 
\usepackage{amsmath} 
\usepackage{amssymb}  
\usepackage{amsfonts}  
\usepackage{epstopdf}
\usepackage{enumerate}
\usepackage{graphicx} 
\usepackage{comment} 
\usepackage{apptools}
\usepackage{subcaption}
\usepackage{theoremref}
\usepackage{bbm}

\usepackage{soul}%
\usepackage{cancel}

\newcommand{\tcb}{\textcolor{black}}

\newtheorem{thm}{Theorem}
\newtheorem{prop}{Proposition}
\newtheorem{lemma}{Lemma}
\newtheorem{corollary}{Corollary}
\newtheorem{definition}{Definition}
\newtheorem{asmp}{Assumption}

\newtheorem{example}{Example}
\newtheorem{remark}{Remark}
\graphicspath{{./},{./figs/}}

\newcommand*{\QEDB}{\hfill\ensuremath{\square}}
    
\begin{document}

\title{
Conditioning with Conditionals
}

\title{On Lie-Bracket Averaging for a Class of Hybrid Dynamical Systems with Applications to Model-Free Control and Optimization}

\author{Mahmoud Abdelgalil and Jorge I. Poveda
\thanks{This work was supported in part by the grants NSF ECCS CAREER 2305756 and AFOSR YIP: FA9550-22-1-0211.}
\thanks{M. Abdelgalil and J. I. Poveda are with the ECE Department, UCSD, CA, USA. Corresponding Author: J. I. Poveda ({\tt poveda@ucsd.edu}).}}

\maketitle

\begin{abstract}
The stability of dynamical systems with oscillatory behaviors and well-defined average vector fields has traditionally been studied using averaging theory. These tools have also been applied to hybrid dynamical systems, which combine continuous and discrete dynamics. However, most averaging results for hybrid systems are limited to first-order methods, hindering their use in systems and algorithms that require high-order averaging techniques, such as hybrid Lie-bracket-based extremum seeking algorithms and hybrid vibrational controllers. To address this limitation, we introduce a novel high-order averaging theorem for analyzing the stability of hybrid dynamical systems with high-frequency periodic flow maps. These systems incorporate set-valued flow maps and jump maps, effectively modeling well-posed differential and difference inclusions. By imposing appropriate regularity conditions, we establish results on $(T,\varepsilon)$-closeness of solutions and semi-global practical asymptotic stability for sets. These theoretical results are then applied to the study of three distinct applications in the context of hybrid model-free control and optimization via Lie-bracket averaging. 
\end{abstract}
%
\section{Introduction}
The theory of averaging has been crucial over the past century for analyzing and synthesizing systems with fast oscillatory or time-varying dynamics \cite{sanders2007averaging}, including nonlinear controllers \cite{khalil}, estimation methods, and model-free optimization algorithms \cite{KrsticBookESC,TanAndNesic2006Local,PoTe17Auto,scheinker2017model}. Stability analyses based on averaging theory have been extensively explored for systems modeled as Ordinary Differential Equations (ODEs) \cite{sanders2007averaging,maggia2020higher} and extended to certain Hybrid Dynamical Systems (HDS) \cite{Wang:12_Automatica,TeelNesicAveraging,PovedaNaliAuto20}, 
 {as well systems wherein dithers are injected into a non-smooth/switched system, and averaging is used to obtain a  continuous approximation of the dynamics \cite{iannelli2003dither,iannelli2006averaging}.}

Stability results based on averaging theory typically rely on a well-defined ``average system'' whose solutions remain ``close" to those of the original system over compact time domains and state space subsets. This closeness is established by transforming the original dynamics into a perturbed version of the average dynamics. By leveraging uniform stability properties via $\mathcal{K}\mathcal{L}$ bounds, it can be shown that the trajectories of the original dynamics converge similarly to those of the average dynamics as the time scale separation increases. {This method has been applied to dynamical systems modeled as ODEs \cite{sanders2007averaging,khalil,TeelNesicAveraging} and some classes of HDS and switching systems \cite{Wang:12_Automatica,liberzon2022stability,PovedaNaliAuto20,de2018hybrid,iannelli2008subtleties,ochoa2023momentum}.} \tcb{These findings have led to the development of new hybrid algorithms for extremum-seeking (ES) and adaptive systems \cite{PoTe17Auto,Kutadinata_Moase,krilavsevic2023learning,abdelgalil2023multi,PatentBenosmanPoveda}.} 

Despite significant progress over the last decade, gaps remain between the averaging tools available for ODEs and those suitable for HDS. Most averaging results for HDS have been limited to ``first-order" averaging, where the nominal stabilizing vector field is derived by neglecting higher-order perturbations. As indicated in \cite[Section 2.9]{sanders2007averaging}, \cite{abdelgalil2022recursive}, and \cite{DurrLieBracket}, first-order approximations may not accurately characterize the stability of certain highly oscillatory systems, such as those in vibrational control \cite{scheinker2017model} and Lie-bracket ES \cite{DurrLieBracket,grushkovskaya2018class}. In these cases, the stabilizing average dynamics require consideration of higher-order terms. This high-order averaging, well-known in ODE literature, is the focus of this paper, which aims to develop similar tools for hybrid dynamical systems with highly oscillatory flow maps and demonstrate their applications in control and optimization problems.

Based on this, our main contributions are the following:

(a) First, a novel high-order averaging theorem is introduced to study the ``closeness of solutions'' property between certain HDS with periodic flow maps and their corresponding average dynamics.  {In contrast to existing literature, we focus on hybrid systems with differential and difference inclusions, where stabilizing effects are determined by an average hybrid system obtained through the recursive application of averaging at high orders of the inverse frequency.} These higher-order averages become essential when traditional first-order average systems fail to capture stability properties of the system. For compact time domains and compact sets of initial conditions, it is established that each solution of the original HDS is $(T,\varepsilon)$-close to some solution of its average hybrid dynamics.  {This notion of closeness between solutions is needed because, in general, HDS generate solutions with discontinuities, or ``jumps'', for which the standard (uniform) distance between solutions is not a suitable metric, see \cite{TeelNesicAveraging,Wang:12_Automatica} and \cite[Ex. 3, pp. 43]{CSM_hybrid}.}  {Moreover, the models studied in this paper admit multiple solutions from a given initial condition, a property that is particularly useful to study families of functions of interest as in e.g., switching systems. To establish the closeness of solutions property, we introduce a new change of variable that is recursive in nature, and which leverages the structure in which the ``high-order'' oscillations appear in the flow map.}

(b) Next, by leveraging the property of closeness of solutions between the original and the average hybrid dynamics, as well as a uniform asymptotic stability assumption on the (second-order) average hybrid dynamics, we establish a semi-global practical asymptotic stability result for the original hybrid system. For the purpose of generality, we allow the average hybrid dynamics to have a compact set that is locally asymptotically stable with respect to some basin of attraction, which recovers the entire space whenever the stability properties are actually global.  {By using proper indicators in the stability analysis, the original hybrid system is shown to render the same compact set semi-global practically asymptotically stable with respect to the same basin of attraction, thus effectively extending \cite[Thm. 2]{TeelNesicAveraging} to high-order averaging.}

c) Finally, by  exploiting the structure of the HDS, as well as the high-order averaging results, we study three different novel stabilization and optimization problems that  {involve hybrid set-valued dynamics and highly oscillatory control:} (a) distance-based synchronization in networks with switching communication topologies and control directions,  {extending the results of \cite{durr2013examples,DurrManifold} to oscillators on dynamic graphs, operating under intermittent control}; (b) source-seeking problems in non-holonomic vehicles operating under \emph{faulty and spoofed sensors}, which extends the results of \cite{DurrLieBracket} to vehicles operating in adversarial environments, and those of \cite{Poveda:20TAC} to vehicles with non-holonomic models with angular actuation; and (c) the solution of model-free \emph{global} extremum seeking problems on smooth compact manifolds,  {which extends the recent results of \cite{ochoa2025robust} to Lie-bracket-based algorithms}.  {In this manner, an important gap in the hybrid extremum-seeking control literature is also filled by integrating hybrid dynamics into Lie-bracket-based extremum-seeking algorithms.}

The rest of this paper is organized as follows. Section \ref{sec_notation} introduces the preliminaries. Section \ref{secmotivationalexample} presents a motivational example.  {Section \ref{closenesstrajectories} presents the main results.} Applications are studied in Section \ref{sec:switching_systems}.  {The proofs are presented in Section \ref{sec:proofs}}, and Section \ref{seconclusions} presents the conclusions.
\section{Preliminaries}
\label{sec_notation}
\subsection{Notation}
 {The bilinear form $\langle \cdot,\cdot\rangle:\mathbb{R}^n\times\mathbb{R}^n\rightarrow\mathbb{R}$ is the canonical inner product on $\mathbb{R}^n$, i.e. $\langle x, y\rangle = x^\top y$, for any two vectors $x,y\in\mathbb{R}^n$.} Given a compact set $\mathcal{A}\subset\mathbb{R}^n$ and a vector $x\in\mathbb{R}^n$, we use $|x|_{\mathcal{A}}:=\min_{\tilde{x}\in\mathcal{A}}\|x-\tilde{x}\|_2$. \tcb{We also use $\overline{\text{con}}(\mathcal{A})$ to denote the closure of the convex hull of $\mathcal{A}$.} A set-valued mapping $M:\mathbb{R}^p\rightrightarrows\mathbb{R}^n$ is outer semicontinuous (OSC) at $z$ if for each sequence $\{z_i,s_i\}\to(z,s)\in\mathbb{R}^p\times\mathbb{R}^n$ satisfying $s_i\in M(z_i)$ for all $i\in\mathbb{Z}_{\geq0}$, we have $s\in M(z)$. A mapping $M$ is locally bounded (LB) at $z$ if there exists an open neighborhood $N_z\subset\mathbb{R}^p$ of $z$ such that  $M(N_z)$ is bounded. The mapping $M$ is OSC and LB relative to a set $K\subset\mathbb{R}^p$ if $M$ is OSC for all $z\in K$ and $M(K):=\cup_{z\in K}M(x)$ is bounded. A function $\beta:\mathbb{R}_{\geq0}\times\mathbb{R}_{\geq0}\to\mathbb{R}_{\geq0}$ is of class $\mathcal{KL}$ if it is nondecreasing in its first argument, nonincreasing in its second argument, $\lim_{r\to0^+}\beta(r,s)=0$ for each $s\in\mathbb{R}_{\geq0}$, and  $\lim_{s\to\infty}\beta(r,s)=0$ for each $r\in\mathbb{R}_{\geq0}$.  Throughout the paper, for two (or more) vectors $u,v \in \mathbb{R}^{n}$, we write $(u,v)=[u^{\top},v^{\top}]^{\top}$. Let $x=(x_1,x_2,\ldots,x_n)\in\mathbb{R}^n$. Given a Lipschitz continuous function $f:\mathbb{R}^n\to\mathbb{R}^m$, we use $\partial_{x_i}f$ to denote the generalized Jacobian \cite{clarke2008nonsmooth} of $f$ with respect to the variable $x_i$. A function $f$ is said to be of class $\mathcal{C}^k$ if its $k$th-derivative is locally Lipschitz continuous. For two vectors $x_1,x_2 \in\mathbb{R}^n$, the notation $x_1\otimes x_2$ denotes the outer product defined as $x_1\otimes x_2 = x_1x_2^{\top}$. We let $\mathbb{S}^1:=\left\{(x_1,x_2)\in\mathbb{R}^2:\,x_1^2+x_2^2 = 1\right\}$ and \tcb{$\mathbb{B}:=\{x\in\mathbb{R}^{n}~:~\|x\|\leq 1\}$, where the dimension $n$ will be clear from the context.} {Finally, all simulation time units may be assumed to be in seconds.}

\vspace{-0.2cm}
\subsection{Hybrid Dynamical Systems}
\subsubsection{Model}
We consider HDS aligned with the framework of \cite{bookHDS}, given by the following inclusions:
\begin{subequations}\label{eq:HDS0}
\begin{align}
&x\in C,~~~~~\dot{x}\in F(x),\\
&x\in D,~~~x^+\in G(x),
\end{align}
\end{subequations}
where $F:\mathbb{R}^n\rightrightarrows\mathbb{R}^n$ is called the flow map, $G:\mathbb{R}^n\rightrightarrows\mathbb{R}^n$ is called the jump map, $C\subset\mathbb{R}^n$ is called the flow set, and $D\subset\mathbb{R}^n$ is called the jump set. We use $\mathcal{H}=(C,F,D,G)$ to denote the \emph{data} of the HDS $\mathcal{H}$. We note that systems of the form \eqref{eq:HDS0} generalize purely continuous-time systems (obtained when $D=\emptyset$) and purely discrete-time systems (obtained when $C=\emptyset$). Of particular interest to us are time-varying systems, which can also be represented as \eqref{eq:HDS0} by using an auxiliary state $\tau\in\mathbb{R}_{\geq0}$ with dynamics $\dot{\tau}=\rho$ and $\tau^+=\tau$, where $\rho\in\mathbb{R}_{\geq0}$ indicates the rate of change of $\tau$. In this paper, we will always work with well-posed HDS that satisfy the following standing assumption \cite[Assumption 6.5]{bookHDS}.
\begin{asmp}\label{asmp:regularity}
The sets $C,D$ are closed. The set-valued mapping $F$ is OSC, LB, and for each $x\in C$ the set $F(x)$ is convex and nonempty. The set-valued mapping $G$ is OSC, LB, and for each $x\in D$ the set $G(x)$ is nonempty.
\end{asmp}
\subsubsection{Properties of Solutions} Solutions to  \eqref{eq:HDS0} are parameterized by a continuous-time index $t\in\mathbb{R}_{\geq0}$, which increases continuously during flows, and a discrete-time index $j\in\mathbb{Z}_{\geq0}$, which increases by one during jumps. Therefore, solutions to \eqref{eq:HDS0} are defined on \emph{hybrid time domains} (HTDs). A set $E\subset\mathbb{R}_{\geq0}\times\mathbb{Z}_{\geq0}$ is called a \textsl{compact} HTD if $E=\cup_{j=0}^{J-1}([t_j,t_{j+1}],j)$ for some finite sequence of times $0=t_0\leq t_1\ldots\leq t_{J}$. The set $E$ is a HTD if for all $(T,J)\in E$, $E\cap([0,T]\times\{0,\ldots,J\})$ is a compact HTD. 

The following definition formalizes the notion of solution to HDS of the form \eqref{eq:HDS0}.
\begin{definition}\label{definitionsolutions1}
A hybrid arc $x$ is a function defined on a HTD. In particular, $x:\text{dom}(x)\to \mathbb{R}^n$ is such that $x(\cdot, j)$ is locally absolutely continuous for each $j$ such that the interval $I_j:=\{t:(t,j)\in \text{dom}(x)\}$ has a nonempty interior. A hybrid arc $x:\text{dom}(x)\to \mathbb{R}^n$ is a solution $x$ to the HDS \eqref{eq:HDS0} if $x(0, 0)\in C\cup D$, and:
\begin{enumerate}
\item For all $j\in\mathbb{Z}_{\geq0}$ such that $I_j$ has nonempty interior: $x(t,j)\in C$ for all $t\in I_j$, and $\dot{x}(t,j)\in F(x(t,j))$ for almost all $t\in I_j$.
\item For all $(t,j)\in\text{dom}(x)$ such that $(t,j+1)\in \text{dom}(x)$: $x(t,j)\in D$ and $x(t,j+1)\in G(x(t,j))$.
\end{enumerate}
A solution $x$ is \emph{maximal} if it cannot be further extended. A solution $x$ is said to be \emph{complete} if $\text{length}~\text{dom}(x)=\infty$. \QEDB
\end{definition}

In this paper, we also work with an ``inflated'' version of \eqref{eq:HDS0}, which is instrumental for robustness analysis \cite[Def. 6.27]{bookHDS}. 
\begin{definition}\label{def:inflatedHDS}
For $\delta>0$, the $\delta$-inflation $\mathcal{H}_{\delta}$ of the HDS $\mathcal{H}$ with data \eqref{eq:HDS0} is given by $\mathcal{H}_{\delta}:=(C_{\delta},D_{\delta},F_{\delta},G_{\delta})$, where the sets $C_{\delta},D_{\delta}$ are defined as:
\begin{subequations}
\begin{align}
C_{\delta}&:=\left\{x\in\mathbb{R}^n:\left(\{x\}+\delta\mathbb{B}\right)\cap C\neq \emptyset \right\},\label{inflatedflowset}\\
D_{\delta}&:=\left\{x\in\mathbb{R}^n:\left(\{x\}+\delta\mathbb{B}\right)\cap D\neq \emptyset\right\},\label{inflatedjumpset}
\end{align}
and the set-valued mappings $F_{\delta},G_{\delta}$ are defined as:
\begin{align}
F_{\delta}(x)&:=\overline{\text{con}}~F((x+\delta\mathbb{B})\cap C)+\delta\mathbb{B},\label{inflatioflowmap}\\
G_{\delta}(x)&:= \{v\in\mathbb{R}^n:~v\in g+\delta\mathbb{B},g\in G((x+\delta\mathbb{B})\cap D)\},\label{inflatiojumpmap}
\end{align}
\end{subequations}

\vspace{-0.2cm}\noindent 
for all $x\in\mathbb{R}^n$. \QEDB
\end{definition}

\vspace{0.1cm}
The following definition, corresponding to \cite[Def. 5.23]{bookHDS}, will be used in this paper.

\vspace{0.1cm}
\begin{definition}\label{def:tauepsilonclose}
Given $T,\rho>0$, two hybrid arcs $x_1:\text{dom}(x_1)\to\mathbb{R}^n$ and $x_2:\text{dom}(x_2)\to\mathbb{R}^n$ are said to be $(T,\rho)$-close if:
\begin{itemize}
\item  for each $(t, j)\in\text{dom}(x_1)$ with $t+j\leq T$ there exists $s$ such that $(s,j)\in\text{dom}(x_2)$, with $|t-s|\leq\rho$ and $|x_1(t,j)-x_2(s,j)|\leq \rho$.
\item for each $(t,j)\in\text{dom}(x_2)$ with $t+j\leq T$ there exists $s$ such that $(s,j)\in\text{dom}(x_1)$, with $|t-s|\leq\rho$ and $|x_2(t,j)-x_1(s,j)|\leq \rho$. \hfill \QEDB  
\end{itemize}  
\end{definition}

\vspace{0.2cm}
\subsubsection{Stability Notions} To study the stability properties of the HDS \eqref{eq:HDS0}, we make use of the following notion, which is instrumental for the study of ``local'' stability properties.

\vspace{0.1cm}
\begin{definition}\label{defproperindicator}
    Let $\mathcal{A}$ be a compact set contained in an open set $\mathcal{B}_{\mathcal{A}}$. A function $\omega:\mathcal{B}_{\mathcal{A}}\rightarrow\mathbb{R}_{\geq 0}$ is said to be a \emph{proper indicator function} for $\mathcal{A}$ on $\mathcal{B}_{\mathcal{A}}$ if $\omega$ is continuous, $\omega(x)=0$ if and only if $x\in\mathcal{A}$, and if the sequence $\{x_i\}_{i=1}^\infty$, $x_i\in\mathcal{B}_{\mathcal{A}}$, approaches the boundary of $\mathcal{B}_{\mathcal{A}}$ or is unbounded then, the sequence $\{\omega(x_i)\}_{i=1}^\infty$ is also unbounded.  \QEDB 
\end{definition}

The use of proper indicators is common when studying (uniform) local stability properties \cite{sontag2022remarks}. The following definition is borrowed from \cite[Def. 7.10 and Thm. 7.12]{bookHDS}.

\begin{definition}\label{definitionstablity1}
    A compact set $\mathcal{A}\subset\mathbb{R}^n$ contained in an open set $\mathcal{B}_{\mathcal{A}}$ is said to be \textit{uniformly asymptotically stable (UAS)} with a basin of attraction $\mathcal{B}_{\mathcal{A}}$ for the HDS \eqref{eq:HDS0} if for every proper indicator $\omega$ on $\mathcal{B}_{\mathcal{A}}$ there exists $\beta\in\mathcal{KL
    }$ such that for each solution $x$ to \eqref{eq:HDS0} starting in $\mathcal{B}_{\mathcal{A}}$, we have
    \begin{align}\label{KLbound1}
        \omega(x(t,j))&\leq \beta(\omega(x(0,0)),t+j), 
    \end{align}
    for all $(t,j)\in\text{dom}(x)$. If $\mathcal{B}_{\mathcal{A}}$ is the whole space, then $\omega(x)=|x|_{\mathcal{A}}$, and in this case the set $\mathcal{A}$ is said to be \emph{uniformly globally asymptotically stable} (UGAS) for the HDS \eqref{eq:HDS0}. \QEDB 
\end{definition}

\vspace{-0.1cm}
\section{Motivational Example}
\label{secmotivationalexample}
 {The problem of target seeking in mobile robotic systems is ubiquitous across various engineering applications, including rescue missions \cite{azzollini2020extremum}, gas leak detection, and general autonomous vehicles performing exploration missions in hazardous environments that could be too dangerous for humans \cite{zhang2007source,Poveda:20TAC}. Different algorithms have been considered for resolving source-seeking problems in mobile robots \cite{zhang2007source}. Naturally, these algorithms rely on real-time exploration and exploitation mechanisms that are robust with respect to different types of disturbances, including small bounded measurement noise and implementation errors \cite{Poveda:20TAC}, as well as certain structured and bounded additive perturbations \cite{suttner2022robustness}. Nevertheless, in many realistic applications, the seeking robots operate in environments where their sensors rarely have continuous and perfect access to measurements of the environment. Common reasons for faulty measurements include intermittent communication networks, interference due to obstacles or external environmental conditions, to name just a few \cite{labar2022extremum,galarza2021extremum}. In addition to encountering faulty and intermittent measurements, autonomous robots operating in adversarial environments might also be subject to malicious \emph{spoofing} attacks that deliberately modify some of the signals used by the controllers. In such situations, it is natural to ask whether the mobile robots can still complete their missions, and under what conditions (if any) such success can be guaranteed.}

 {To study the above question, we consider a typical model of a mobile robot studied in the context of source seeking problems: a planar kinematic model of a non-holonomic vehicle, with equations
\begin{align*}
    \dot{x}_1&= u\,x_3, & \dot{x}_2&= u\,x_4, & \dot{x}_3&= \Omega x_4, & \dot{x}_4&=-\Omega x_3,
\end{align*}
where the vector $x_p:=(x_1,x_2)\in\mathbb{R}^2$ models the position of the vehicle in the plane, the vector $x_{3,4}:=(x_3,x_4)\in\mathbb{S}^1$ captures the orientation of the vehicle, $u$ is the forward velocity, and $\Omega$ is the angular velocity.  {We use the above model for the kinematics of the vehicle to simplify the discussion since our focus is on modeling the effect of intermittent measurements and spoofing attacks. Nevertheless, this model is ubiquitous in the source-seeking literature, see \cite{Pappas_SourceSeeking}, \cite{CochranSourceSeekingTAC09}, \cite{Poveda:20TAC}, \cite{zhang2007source}, and \cite{durrEbenbauer}.}} 

 {The main objective of the vehicle is to stabilize its position at a particular unknown target point $x_p^{\star}\in\mathbb{R}^2$ using only real-time measurements of a distance-like signal $J(x_p)$, which for simplicity is assumed to be $\mathcal{C}^1$, strongly convex with minimizer at the point $x_p^{\star}$, and having a globally Lipschitz gradient (these assumptions will be relaxed later in Section \ref{intermittentliebracket}). \tcb{In other words, the vehicle seeks to solve the problem $\min_{x_p\in\mathbb{R}^2}~J(x_p)$ without knowledge of the mathematical form of $J$.} Since nonholonomic vehicles cannot be simultaneously stabilized in position and orientation using smooth feedback, we let the vehicles oscillate continuously using a predefined constant angular velocity $\Omega > 0$. This model can also capture vehicles that have no \textit{directional} control over $\Omega$, such as in quad-copters that have lost one propeller, leading to a constant yaw rate during hover that can be effected, though not eliminated or reversed, by modulating the remaining propellers \cite{mueller2014stability}.} 

\begin{figure}[t!]
    \centering    \includegraphics[width=0.475\textwidth]{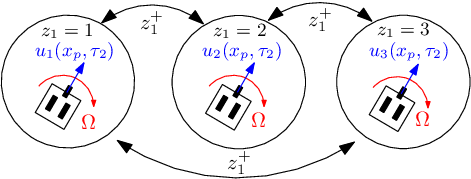}
    \caption{\small{Automaton-like representation of switching seeking dynamics in nonholonomic vehicles with multiple operating modes: under spoofing ($z_1=1$), under no measurement ($z_1=2$), and under nominal operation ($z_1=3$).}}
\label{fig:automaton}
\vspace{-0.3cm}
\end{figure}
Under ``nominal'' operating conditions, and to stabilize (a neighborhood of) the target $x_p^*$, we can consider the following feedback controller that only uses real-time measurements of the potential field $J$ evaluated at the current position
\begin{align*}
    u(x_p,\tau_2) &= \varepsilon^{-1} \cos\left(\tau_2 + J(x_p)\right),~\dot{\tau}_2=\frac{1}{\varepsilon^2},~\Omega=\frac{1}{\varepsilon},
\end{align*}
where $\varepsilon>0$ is a tunable parameter. To capture the effect of having intermittent measurements and malicious spoofing attacks in the algorithm, we re-write the nominal closed-loop system as the following mode-dependent dynamical system:
\begin{subequations}\label{nonholonomicvehiclea1}
 \begin{align}
    \dot{x}_1&= \varepsilon^{-1}x_3u_{z_1}(x_p,\tau_2), & \dot{x}_3&= \varepsilon^{-1}x_4, \\
    \dot{x}_2&= \varepsilon^{-1}x_4u_{z_1}(x_p,\tau_2), & \dot{x}_4&= -\varepsilon^{-1}x_3, \\
    \dot{\tau}_2&= \varepsilon^{-2},& \dot{z}_1&= 0,
\end{align}   
where $z_1\in\mathcal{Q}:=\{1,2,3\}$ is a logic mode that is kept constant during the evolution of \eqref{nonholonomicvehiclea1}. In this case, the map $u_{z_1}$ is
\begin{align}\label{controlawtoy}
    u_{z_1}(x_p,\tau_2)=  \cos\left(\tau_2 + (z_1-2)J(x_p)\right),
\end{align}
\end{subequations} 
where $z_1=3$ describes the nominal operation mode, $z_1=2$ describes the mode in which no measurement of $J$ is available to the vehicle, and $z_1=1$ corresponds to the mode in which the control algorithm is under spoofing that is able to reverse the sign of the measurements of $J$.  {In this way, the dynamics of the vehicle switch in real-time between the three operating modes as the robot simultaneously seeks the target $x_p^*$.} Note that, at every switching instant, the states $x=(x_1,x_2,x_3,x_4)$ remain constant, and $z_1$ jumps according to $z_1^+\in \mathcal{Q}\backslash\{z_1\}$, i.e., the system is allowed to switch from its current mode to any of the other two modes. Figure \ref{fig:automaton} illustrates this switching behavior.
To study the stability properties of the system, we consider the coordinate transformation:
\begin{align*}
    \begin{pmatrix}
        {x}_3\\{x}_4
    \end{pmatrix}&= \begin{pmatrix}
        \cos\left(\varepsilon^{-1}t\right) & \sin\left(\varepsilon^{-1}t\right)\\
        -\sin\left(\varepsilon^{-1}t\right) & \cos\left(\varepsilon^{-1}t\right)
    \end{pmatrix}\begin{pmatrix}
        \tilde{x}_3\\\tilde{x}_4
    \end{pmatrix},
\end{align*}
which is defined via a rotation matrix and therefore is invertible for all $t\geq0$. Using $\tau_1=\varepsilon^{-1}t$, the transformation leads to the new system
\begin{subequations}\label{eq:mot_example_org_sys_shifted}
\begin{align}
    \dot{x}_p &=\varepsilon^{-1} \tcb{\underbrace{\begin{pmatrix}
        \cos\left(\tau_1\right) & \sin\left(\tau_1\right)\\
        -\sin\left(\tau_1\right) & \cos\left(\tau_1\right)
    \end{pmatrix}}_{A(\tau_1)}}\begin{pmatrix}
        \tilde{x}_3\\\tilde{x}_4
    \end{pmatrix}u_{z_1}(x_p,\tau_2),\\
    \dot{\tilde{x}}_{3,4}&=0, \quad
    \dot{\tau}_1=\varepsilon^{-1}, \quad \dot{\tau}_2=\varepsilon^{-2},~~~\dot{z}_1=0,
\end{align}
\end{subequations}
with jumps $x^+=x$ and $z_1^+\in \mathcal{Q}\backslash\{z_1\}$. In \eqref{eq:mot_example_org_sys_shifted}, the vector $\tilde{x}_{3,4}$ is now constant, but still restricted to the compact set $\mathbb{S}^1$. For this system, we can investigate the stability properties of the state $\tilde{x}=(x_p,\tilde{x}_{3,4})$ with respect to the set $\{x_p^\star\}\times\mathbb{S}^1$. To do this, we model the closed-loop system as a HDS of the form \eqref{eq:HDS0}, with switching signal $z_1$ being generated by an extended auxiliary hybrid automaton, with states $z:=(z_1,z_2,z_3)\in\mathbb{R}_{\geq0}^3$, and set-valued dynamics 
\begin{subequations}\label{hybridautomaton}
\begin{align}
&z\in C_{z},~~\dot{z}\in \mathcal{T}_F(z):=\left(\begin{array}{c}
\{0\}\\
\left[0,\eta_1\right]\\
\left[0,\eta_2\right]-\mathbb{I}_{\mathcal{Q}_u}(z_1)
\end{array}\right),\label{flowsautomatonswitching}\\
&z\in D_{z},~~z^+\in \mathcal{T}_G(z):=\left(\begin{array}{c}
\mathcal{Q}\backslash\{z_1\}\\
z_2-1\\
z_3\\
\end{array}\right),
\end{align}
\end{subequations}
where $\eta_1>0$,  {$\eta_2\in[0,1)$}, $\mathbb{I}_{\mathcal{Q}_u}(\cdot)$ is the classic indicator function, the set of logic modes has been partitioned as $\mathcal{Q}=\mathcal{Q}_s\cup\mathcal{Q}_u$, with $\mathcal{Q}_s=\{3\}$ and $\mathcal{Q}_u=\{1,2\}$, and the sets $C_{z},D_{z}$ are given by
\begin{equation}\label{setsauxiliary001}
C_{z}=\mathcal{Q}\times [0,N_\circ]\times[0,T_\circ],~~D_{z}=\mathcal{Q}\times [1,N_\circ]\times[0,T_\circ].
\end{equation}
\tcb{for some constants $N_\circ\geq 1$ and $T_\circ>0$. In this system, the states $z_2$ and $z_3$ can be seen as timers used to coordinate when the system jumps.}
\begin{figure}
    \centering    \includegraphics[width=0.475\textwidth]{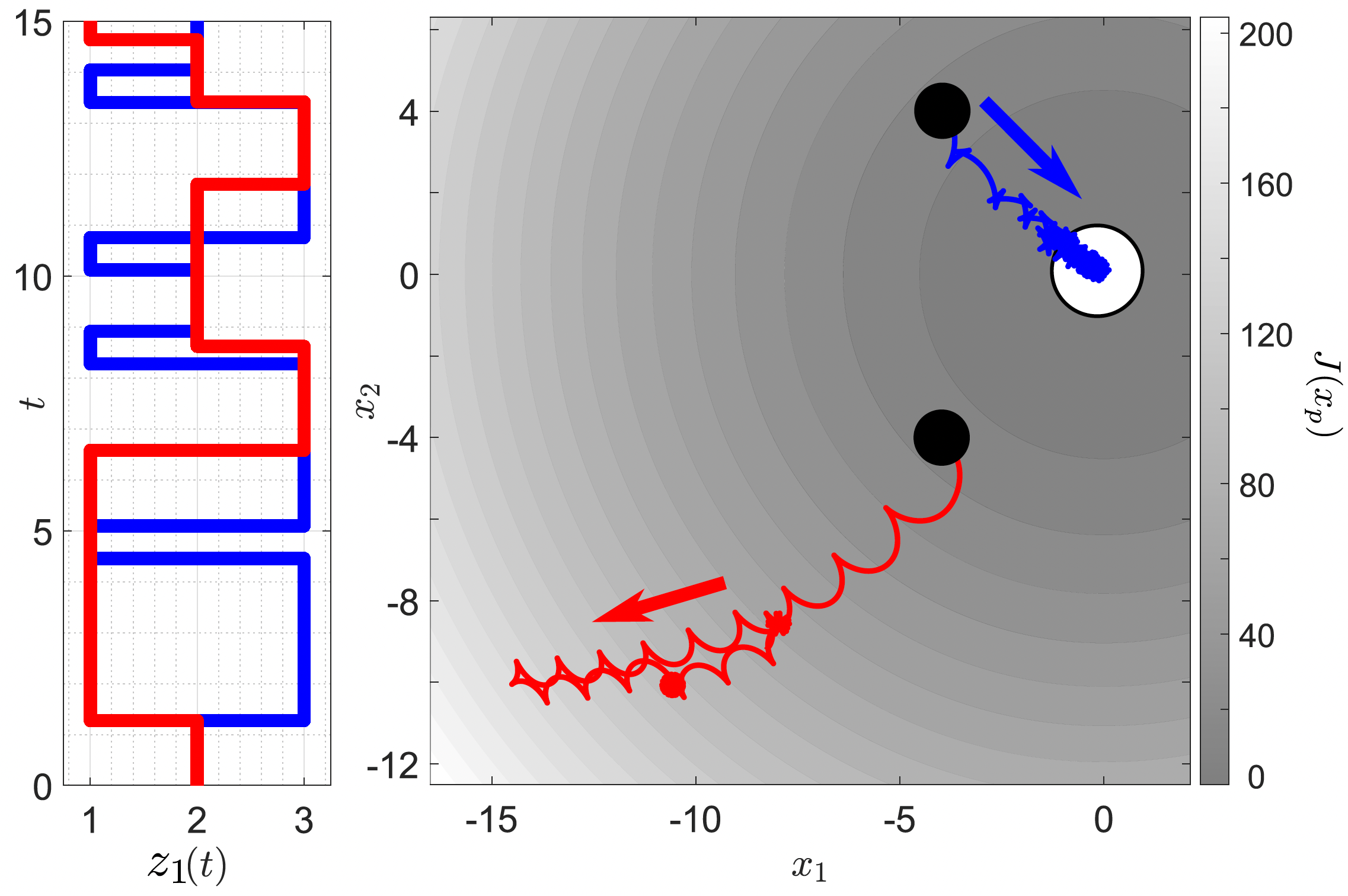}
    \caption{\small{ {Trajectories of the vehicle starting from two initial conditions $(x_1(0),x_2(0))=(-4,4)$ (blue) and $(x_1(0),x_2(0))=(-4,-4)$ (red), for $\varepsilon = 1/\sqrt{10\pi}$, $J(x_p)=0.5 x_1^2+0.5x_2^2$, and for two different switching signals $z_1(t)$, shown in the left plot. The black dots indicate the initial conditions.}}}
    \label{fig:source_seeking_example_numerical_simulation}
    \vspace{-0.3cm}
\end{figure}
 {In particular, every hybrid arc generated by the hybrid automaton \eqref{hybridautomaton} satisfies the following two properties for any two times $t_2>t_1$ that are in its domain \cite[Ex. 2.15]{bookHDS},\cite[Lemma 7]{PoTe17Auto}:
\begin{subequations}\label{switchingconditionsineq}
\begin{align}
    N_\sharp(t_1,t_2)&\leq \eta_1(t_2-t_1) + N_\circ, \label{ADT}\\
    T_\sharp(t_1,t_2)&\leq \eta_2(t_2-t_1) + T_\circ,\label{ATT}
\end{align}
\end{subequations}
where $N_\sharp(t_1,t_2)$ is the total number of jumps during the time interval $(t_1,t_2)$, and $T_\sharp(t_1,t_2)=\int_{t_1}^{t_2}\mathbb{I}_{\mathcal{Q}_u}(z_1(t))\,dt$ is the total activation time during $[t_1,t_2]$, and during flows of the system, of the modes in the set $\mathcal{Q}_u$. \tcb{In fact, condition \eqref{ADT} imposes an average dwell-time (ADT) constraint \cite{AverageDwellTime} on the switches of $z_1$, while condition \eqref{ATT} imposes an average activation time (AAT) constraint \cite{GuosongLiberzon} on the time spent in modes 1 and 2.} These conditions are parameterized by the tunable constants $\eta_1$ and $\eta_2$, respectively. Note that \eqref{ADT} immediately rules out Zeno behavior.}  {By modeling the switching signals as solutions of the hybrid automaton \eqref{hybridautomaton}, we can study the closed-loop system without pre-specifying the switching times of $z_1$,  {which are generally unknown.} Instead, we consider any possible solution $(\tilde{x},z)$ generated by the interconnection \eqref{eq:mot_example_org_sys_shifted}-\eqref{setsauxiliary001}.}

 {While the $x_p$-dynamics in \eqref{eq:mot_example_org_sys_shifted} are periodic with high-frequency oscillations, existing averaging tools in the literature \cite{TeelNesicAveraging} do not capture the stabilizing effect of the control law \eqref{controlawtoy}. This can be observed by introducing a new time scale $s=\frac{t}{\varepsilon}$, which leads to the following dynamics (in the s-time scale):
\begin{equation}\label{averagebad}
\dot{x}_p=A(\tau_1)\tilde{x}_{3,4}u_{z_1}(x_p,\tau_2),~~\dot{\tau}_1=1,~~\dot{\tau}_2=\frac{1}{\varepsilon}.
\end{equation}
Using \eqref{controlawtoy} and the property $\cos(\alpha+\beta)=\cos(\alpha)\cos(\beta)-\sin(\alpha)\sin(\beta)$, it is easy to see that for any constant vector $\tilde{x}_{3,4}\in\mathbb{S}^1$ the average of the vector field \eqref{averagebad} along the fast varying state $\tau_2$ is equal to zero. \tcb{Since the first-order average vector field is zero, which is only marginally stable, no stability conclusions can be obtained from a direct application of first-order averaging theory. Indeed, in this case, the stabilizing effects are dictated by higher-order terms that are neglected from the first-order average.} Similar obstacles emerge when using first-order averaging theory in the analysis of Lie-bracket-based extremum seeking controllers \cite{DurrLieBracket} and in certain vibrational controllers \cite{scheinker2017model}.} On the other hand, as we will show in the next section, by using \emph{second-order} averaging theory for HDS, we can obtain the following second-order average hybrid dynamics of  \eqref{eq:mot_example_org_sys_shifted}-\eqref{setsauxiliary001}, with states $(\bar{x},\bar{z})$:
\begin{align*}
&(\bar{x},\bar{z})\in \left(\mathbb{R}^2\times\mathbb{S}^1\right) \times C_{z},~\left\{\begin{array}{l}
~~\dot{\bar{x}}_p=\dfrac{(2-\bar{z}_1)}{2}\nabla J(\bar{x}_p),\\
\dot{\bar{x}}_{3,4}=0\vspace{0.1cm},\\
~~~~\dot{\bar{z}}\in \mathcal{T}_{F}(\bar{z}),
\end{array}\right.\\
&(\bar{x},\bar{z})\in \left(\mathbb{R}^2\times\mathbb{S}^1\right)\times D_{z},~\left\{\begin{array}{l}
\bar{x}^+=\bar{x},\\
\bar{z}^+\in \mathcal{T}_G(\bar{z}).
\end{array}\right.
\end{align*}
In turn, under the smoothness and strong convexity assumption on $J$, this HDS renders UGAS the compact set $\mathcal{A}:=\{x_p^{\star}\}\times\mathbb{S}^1\times\mathcal{Q}\times [0,N_{\circ}]\times[0,T_{\circ}]$ for $\eta_1>0$ and $\eta_2$ sufficiently small (see Section \ref{sec:switching_systems}-A-2). By using Theorem 2 in Section \ref{closenesstrajectories}, we will be able to conclude that, for any compact set $K_0\subset\mathbb{R}^2$ and any $\nu>0$, there exists $\varepsilon^*>0$ such that for all $\varepsilon\in(0,\varepsilon^*)$ every trajectory of \eqref{eq:mot_example_org_sys_shifted} that starts in $K_0$, satisfies
\begin{equation*}
|x_p(t,j)-x^{\star}|\leq \beta\left(|x_p(0,0)-x^{\star}|,t+j\right)+\nu,
\end{equation*}
for all $(t,j)$ in the domain of the solution. A more general result will be stated later in Section \ref{subsectionswitching} for a class of hybrid Lie-bracket averaging-based systems. In the meantime, we illustrate in Figure \ref{fig:source_seeking_example_numerical_simulation} (in blue color) the convergence properties of the vehicle under the controller \eqref{controlawtoy}. The inset also shows the evolution in time of the switching signal $z_1$, indicating which mode is active at each time. As observed, the vehicle operating under spoofing and intermittent sensing  {($T_\sharp(0,15)\approx 5.98$)} is able to successfully complete the source-seeking mission. We also show in red color an unstable trajectory of the vehicle under the same controller and a more frequent spoofing attack  {($T_\sharp(0,15)\approx 11.34$)} that does not satisfy \eqref{ATT}  {(i.e., with a larger value of $\eta_2$ in \eqref{ATT})}, also shown in the inset. For both scenarios depicted in Figure \ref{fig:source_seeking_example_numerical_simulation}, the switching signals $z_1$ are generated as solutions of the hybrid automaton \eqref{hybridautomaton}.
\begin{remark}
 {The results and tools discussed in this section (extended in Section \ref{subsectionswitching}) can help practitioners evaluate the resilience of non-holonomic vehicle seeking dynamics against intermittent measurements and spoofing attacks. By characterizing the parameters $(\eta_1, \eta_2)$ and $\varepsilon$ that ensure stable (practical) source seeking, practitioners can design effective detection and rejection mechanisms to ensure the control system operates in its nominal mode \tcb{``sufficiently often".}} \QEDB 
\end{remark}

\vspace{0.1cm}
The motivational problem in this section was modeled using a high-amplitude, high-frequency oscillatory system, whose averaged system corresponds to a switching system. However, the results presented in this paper are applicable to a broader class of hybrid systems, with switching systems representing just one specific case.
\section{Second-Order Averaging for a Class of Hybrid Dynamical Systems}
\label{closenesstrajectories}
We consider a subclass of HDS \eqref{eq:HDS0}, given by
\begin{subequations}\label{eq:ms_orig_hybrid_sys}
\begin{align}
\left(({x},z),\tau\right)&\in C\times\mathbb{R}_{\geq0}^2,  & &\begin{cases}    
\left(\begin{array}{c}\dot{{x}}\\\dot{z}\end{array}\right)&\in F_{\varepsilon}({x},{z},\tau_1,\tau_2),\\
    \hphantom{\Big[}\begin{array}{c}\dot{\tau}_1\end{array}\hphantom{\Big]}&= \varepsilon^{-1},\\
    \hphantom{\Big[}\begin{array}{c}\dot{\tau}_2\end{array}\hphantom{\Big]}&= \varepsilon^{-2},
\end{cases}\label{secondcaseHDS0}\\
\left(({x},{z}),\tau\right)&\in D\times\mathbb{R}_{\geq0}^2, & &\begin{cases}    
\left(\begin{array}{c}{x}^+\\{z}^+\end{array}\right)&\in G({x},{z}),\\
\hphantom{\Big[}\begin{array}{c}\tau_1^+\end{array}\hphantom{\Big]} &= \tau_1,\\
\hphantom{\Big[}\begin{array}{c}\tau_2^+\end{array}\hphantom{\Big]} &= \tau_2,
\end{cases}\label{secondcaseHDS}
\end{align}
\end{subequations}
where $\varepsilon>0$ is a small parameter, $x\in\mathbb{R}^{n_1}$, $z\in\mathbb{R}^{n_2}$, $n_1+n_2=n$, $C,D\subset\mathbb{R}^n$ $\tau_i\in\mathbb{R}_{\geq0}$, for $i\in\{1,2\}$, and $\tau=(\tau_1,\tau_2)\in\mathbb{R}_{\geq0}^2$. The set-valued mapping $F_{\varepsilon}$ is given by
\begin{align}\label{flowmapstructure}
F_{\varepsilon}({x},{z},\tau_1,\tau_2):=\big\{f_\varepsilon(x,z,\tau_1,\tau_2)\big\}\times \Phi({z}),
\end{align}
where $f_\varepsilon$ is a continuous function of the form
\begin{align}\label{vectofield0}
f_\varepsilon(x,z,\tau_1,\tau_2)= \sum_{k=1}^2\varepsilon^{k-2}{\phi}_{k}({x},{z},\tau_1,\tau_2).
\end{align}
The functions $\{\phi_k\}_{k=1}^2$, the sets $C,D\subset\mathbb{R}^n$, and the set-valued mappings $\Phi$, $G$ are application-dependent. Their regularity properties are characterized by the following assumption:

\vspace{0.1cm}
\begin{asmp}\label{asmp:basic}
$C$ and $D$ are closed. $\Phi$ is OSC, LB, and convex-valued relative to $C$; $G$ is OSC and LB relative to $D$; for all $z\in C$ the set $\Phi(z)$ is nonempty; and for all $(x,z)\in D$ the set $G(x,z)$ is nonempty.\QEDB
\end{asmp}

\vspace{0.1cm}
For the function $f_{\varepsilon}$ defined in \eqref{vectofield0} we will ask for stronger regularity conditions in terms of smoothness and periodicity of the mappings $\phi_k$ in a small inflation of $C$ (c.f. Eq. \eqref{inflatedflowset}). 
\vspace{0.1cm}
\begin{asmp}\thlabel{asmp:A1}\normalfont
There exists $\delta^*>0$ and $T_k>0$, for $k\in\{1,2\}$, such that for all $\delta\in[0,\delta^*]$ the following holds: 
\begin{enumerate}[(a)]
\item \emph{Smoothness}: For $k\in\{1,2\}$, the functions $\phi_k:C_{\delta}\times\mathbb{R}^2_+\to\mathbb{R}^{n_1}$ satisfy that: $\phi_k\in\mathcal{C}^{2-k}$ with respect to $(x,\tau_1)$, $\phi_k\in\mathcal{C}^{0}$ with respect to $z$, and $\phi_k$ is continuous in $\tau_2$.
\item \emph{Periodicity}: For $k\in\{1,2\}$, $\phi_k$ satisfies
\begin{align*}
{\phi}_k({x},{z},\tau_1+T_1,\tau_2) &= {\phi}_k({x},{z},\tau_1,\tau_2),\\
    {\phi}_k({x},{z},\tau_1,\tau_2+T_2) &= {\phi}_k({x},{z},\tau_1,\tau_2),
\end{align*}
for all $((x,z),\tau)\in C_{\delta}\times\mathbb{R}^2_+$.

\item \emph{Zero-Average of $\phi_1$ in} $\tau_2$: 
\begin{equation*}
{\textstyle\int_0^{T_2}}{\phi}_1({x},{z},\tau_1,s_2)\,ds_2 = 0,
\end{equation*}
for all $((x,z),\tau)\in C_{\delta}\times\mathbb{R}^2_+$. \QEDB 
\end{enumerate}
\end{asmp}

\vspace{0.15cm}
\begin{remark}
\tcb{The periodicity and smoothness assumptions on $\phi_k$ are standard \cite[Ch.~2]{sanders2007averaging}. Periodicity is particularly important in higher-order averaging as it simplifies the computations of the near-identity transformations required to establish the closeness-of-solutions property \cite[Section~2.9]{sanders2007averaging}. Without the periodicity assumption, rigorous proofs become significantly more complex. For example, in the general time-varying case, higher-order averaging necessitates the use of general-order functions of the small parameter rather than polynomials. This increases the complexity of the asymptotic approximations that can be achieved \cite[Theorem~4.5.4]{sanders2007averaging} without providing additional practical benefits, since in the applications of interest, the periodicity assumption holds by design.}\QEDB 
\end{remark}

To study the stability properties of the HDS \eqref{eq:ms_orig_hybrid_sys}, we first characterize its \emph{2$^{\text{nd}}$-order average HDS}. To define this system,  for each $((x,z),\tau)\in C_{\delta}\times\mathbb{R}^2_+$ we define the following auxiliary function: 
\begin{subequations}\label{averagemappingsdefined0}
\begin{align}\label{eq:u1def}
    u_1(x,z,\tau):= \int_0^{\tau_2} \phi_1(x,z,\tau_1,s_2)\,ds_2,
\end{align}
as well as the \emph{Lie-bracket} between \eqref{eq:u1def} and $\phi_1$:
\begin{equation}\label{liebracket}
\left[u_1,\phi_1\right]_x(\cdot)=\left(\partial_x \phi_1 \cdot u_1 - \partial_x u_1 \cdot \phi_1\right)(\cdot),
\end{equation}
\end{subequations}
where we omitted the function arguments to simplify notation. Using \eqref{averagemappingsdefined0}, we introduce the \emph{2$^{\text{nd}}$-order average mapping} of $f_{\varepsilon}$, denoted $\bar{f}:C_{\delta}\to\mathbb{R}^{n_1}$, given by:
\begin{align}\label{fbareq}
\bar{f}(\theta):=&\frac{1}{T_1T_2}\int_0^{T_1}\hspace{-0.25cm}\int_0^{T_2}\hspace{-0.1cm}\left(\phi_2(\theta,\tau)+\frac{1}{2}\left[u_1,\phi_1\right]_x(\theta,\tau)\right)\,d\tau_2\,d\tau_1,
\end{align}
where $\theta:=(x,z)$. The following definition leverages $\bar{f}(\cdot)$.

\vspace{0.1cm}
\begin{definition}\label{defaveragehybridsystem}
The 2$^\text{nd}$-order average HDS of \eqref{eq:ms_orig_hybrid_sys} has states $\bar{\theta}:=(\bar{x},\bar{z})\in\mathbb{R}^{n_1}\times\mathbb{R}^{n_2}$, and is given by 
\begin{equation}\label{eq:ms_avg_hybrid_sys}
\mathcal{H}_2^{\text{ave}}:~~~\left\{\begin{array}{ll}
\bar{\theta}\in C, &   
~~\dot{\bar{\theta}}\in F\left(\bar{\theta}\right)
\vspace{0.2cm}\\
\bar{\theta}\in D, &     
\bar{\theta}^+\in G\left(\bar{\theta}\right)
\end{array}\right.,
\end{equation}
where 
\begin{equation*}
F(\bar{\theta}):=\left\{\bar{{f}}\left(\bar{\theta}\right)\right\}\times \Phi(z),~~~~\forall~\bar{\theta}\in C,
\end{equation*}
and  {$\Phi,C,G,D$ are the same of \eqref{eq:ms_orig_hybrid_sys}.}
\QEDB
\end{definition}

\vspace{0.1cm}

\begin{remark}
It is worth mentioning that if $\phi_1$ satisfies
\vspace{-0.2cm} 
\begin{align*}
    \phi_1(\theta,\tau)&= \sum_{\ell=1}^r{b}_{\ell}({\theta},\tau_1)\,v_\ell(\tau), \\
    \textstyle \int_0^{\tau_2}\phi_1({\theta},\tau_1,s_2)ds_2&= \sum_{\ell=1}^r{b}_{\ell}(\theta,\tau_1)\,\int_0^{\tau_2} v_\ell(\tau_1,s_2)\,ds_2,
\end{align*}
for some functions $b_{\ell},v_{\ell}$ and some $r\in\mathbb{Z}_{\geq1}$, then the Lie bracket in \eqref{liebracket} reduces to:
\begin{align*}
\left[u_1,\phi_1\right]_{x}(\cdot)&=\hspace{-0.2cm}\sum_{\ell_1,\ell_2=1}^r\left[{b}_{\ell_1},{b}_{\ell_2}\right]_{x}(\theta,\tau_1)\int_0^{s_2}\hspace{-0.1cm}v_{\ell_1}(\tau_1,s_2)v_{\ell_2}(\tau)\,ds_2\\
&=\hspace{-0.2cm}\sum_{\ell_1>\ell_2}^r\left[{b}_{\ell_1},{b}_{\ell_2}\right](\theta,\tau_1)v_{\ell}(\tau),
\end{align*}
where $v_{\ell}$ is the time-varying function:
\begin{align*}
v_{\ell}(\tau)=\hspace{-0.1cm}\int_0^{\tau_2}\hspace{-0.1cm}v_{\ell_1}(\tau_1,s_2)v_{\ell_2}(\tau)\,ds_2-\int_0^{\tau_2}\hspace{-0.1cm}v_{\ell_2}(\tau_1,s_2)v_{\ell_1}(\tau)\,ds_2,
\end{align*}
and, consequently, the average of $\left[u_1,\phi_1\right]_x(\cdot)$ reduces to:
\begin{equation*}
\frac{1}{T_2}\int_{0}^{T_2}\left[u_1,\phi_1\right]_x(\theta,\tau_1,s_2)ds_2 = \hspace{-0.2cm}\sum_{\ell_1>\ell_2}^r\left[{b}_{\ell_1},{b}_{\ell_2}\right]_x(\theta,s_1)\Lambda(\tau_1),
\end{equation*}
where $\Lambda(\tau_1)=\frac{1}{T_2}\int_{0}^{T_2}v_{\ell}(\tau_1,s_2)ds_2$.
 {However, in contrast to the standard Lie-bracket averaging framework \cite{DurrLieBracket}, here we do not necessarily assume that ${\phi}_1$ admits such a decomposition, which enables the study of more general vector fields $f_{\varepsilon}$.} \QEDB 
\end{remark}

\vspace{0.05cm}
The next lemma follows directly by Assumption \ref{asmp:basic} and the continuity properties of $\phi_k$, $k\in\{1,2\}$, c.f., Assumption \ref{asmp:A1}.

\vspace{0.05cm}
\begin{lemma}
Suppose that Assumptions \ref{asmp:basic}-\ref{asmp:A1} hold. Then, the HDS \eqref{eq:ms_orig_hybrid_sys} and the HDS $\mathcal{H}^{\text{ave}}_2$ in \eqref{eq:ms_avg_hybrid_sys} satisfy Assumption \ref{asmp:regularity}. \QEDB 
\end{lemma}

\vspace{0.05cm}
 {The following Theorem the first main result of this paper. It uses the notion of $(T,\rho)$-closeness, introduced in Definition \ref{def:tauepsilonclose}.} All the proofs are deferred to Section \ref{sec:proofs}.

\vspace{0.05cm}
\begin{thm}\thlabel{thm:closeness_of_trajs}\normalfont
Suppose that Assumptions \ref{asmp:basic}-\ref{asmp:A1} hold. Let ${K}_0\subset\mathbb{R}^n$ be a compact set such that every solution of the HDS \eqref{eq:ms_avg_hybrid_sys} with initial conditions in ${K}_0$ has no finite escape times. Then, for each $\rho>0$ and each $T>0$, there exists $\varepsilon^*>0$ such that for each $\varepsilon\in(0,\varepsilon^*)$ every solution of  (\ref{eq:ms_orig_hybrid_sys}) starting in ${K}_0+\rho\mathbb{B}$ is $(T,\rho)$-close to some solution of (\ref{eq:ms_avg_hybrid_sys}) starting in ${K}_0$. \QEDB 
\end{thm}

\vspace{0.1cm}
\begin{remark}
 {The result of \thref{thm:closeness_of_trajs} establishes a novel closeness-of-solutions property (in the sense of Definition \ref{def:tauepsilonclose}) between the high-frequency high-amplitude periodic-in-the-flows HDS \eqref{eq:ms_orig_hybrid_sys} and the second-order average HDS  $\mathcal{H}_2^{\text{ave}}$ defined in \eqref{eq:ms_avg_hybrid_sys}. Unlike Lie-bracket averaging results for Lipschitz ODEs \cite[Thm. 1]{DurrLieBracket}, the result does not assert closeness between all solutions of \eqref{eq:ms_orig_hybrid_sys} and $\mathcal{H}_2^{\text{ave}}$, but rather that for every solution of the original HDS \eqref{eq:ms_orig_hybrid_sys} there exists a solution of the 2nd-order average dynamics $\mathcal{H}_2^{\text{ave}}$ that is $(T,\rho)$-close. Thus, Theorem \ref{thm:closeness_of_trajs} effectively extends to the 2nd-order case existing 1st-order averaging results for HDS, such as \cite[Thm. 1]{TeelNesicAveraging}.} \QEDB 
\end{remark}

\vspace{0.05cm}
As a result of Theorem \ref{thm:closeness_of_trajs}, if all trajectories \(\bar{\theta}\) of \(\mathcal{H}_2^{\text{ave}}\) satisfy suitable convergence bounds, then the trajectories \(\theta\) of the original HDS \eqref{eq:ms_orig_hybrid_sys} will approximately inherit the same bounds on compact sets and time domains. To extend these bounds to potentially unbounded hybrid time domains, we use the following stability definition \cite[Def. 6]{TeelNesicAveraging}.

\vspace{0.08cm}
\begin{definition}\thlabel{defn:SPUAS}
    For the HDS (\ref{eq:ms_orig_hybrid_sys}), a compact set $\mathcal{A}\subset\mathbb{R}^{n}$ is said to be \emph{Semi-Globally Practically Asymptotically Stable} (SGPAS) with respect to the basin of attraction $\mathcal{B}_{\mathcal{A}}$ as $\varepsilon\rightarrow 0^+$ if, for each proper indicator function $\omega$ on $\mathcal{B}_{\mathcal{A}}$ there exists $\beta\in\mathcal{KL}$ such that, for each compact set $K\subset\mathcal{B}_{\mathcal{A}}$ and for each $\nu>0$, there exists $\varepsilon^*>0$ such that for all $\varepsilon\in(0,\varepsilon^*]$ and for all solutions of \eqref{eq:ms_orig_hybrid_sys} with $\theta(0,0)\in K$ we have
    \begin{align}\label{KLboundSGPAS}
        \omega(\theta(t,j))\leq \beta(\omega(\theta(0,0)),t+j)+ \nu, 
    \end{align}
    for all $(t,j)\in\text{dom}(\theta,\tau)$. \QEDB 
\end{definition}

As discussed in Definition \ref{definitionstablity1}, if $\mathcal{B}_{\mathcal{A}}$ covers the whole space, then we can also take $\omega(x)=|x|_{\mathcal{A}}$ in \eqref{KLboundSGPAS}. If, additionally, $C\cup D$ is compact, then \eqref{KLboundSGPAS} describes a Global Practical Asymptotic Stability (GPAS) property.

 {With \thref{defn:SPUAS} at hand, we can now state the second main result of this paper, which establishes a novel second-order averaging result for HDS:}

\begin{thm}\thlabel{thm:avg_UAS_implies_org_PUAS}
Suppose that Assumptions \ref{asmp:basic} and \ref{asmp:A1} hold, and let $\mathcal{A}\subset\mathbb{R}^n$ be a compact set. If $\mathcal{A}$ is UAS with basin of attraction $\mathcal{B}_{\mathcal{A}}$ for $\mathcal{H}_2^{\text{ave}}$, then $\mathcal{A}$ is also SGPAS as $\varepsilon\to0^+$ with respect to the basin of attraction $\mathcal{B}_{\mathcal{A}}$ for the HDS \eqref{eq:ms_orig_hybrid_sys}. \QEDB 
\end{thm}

The discussion preceding Theorem \ref{thm:avg_UAS_implies_org_PUAS} directly implies the following Corollary, which we will leverage in Section \ref{sec:switching_systems} for the study of \emph{globally} practically stable Lie-bracket hybrid ES systems on smooth boundaryless compact manifolds.

\begin{corollary}\thlabel{corollary1}
Suppose that Assumptions \ref{asmp:basic} and \ref{asmp:A1} hold, and let $\mathcal{A}\subset\mathbb{R}^n$ be a compact set. If $C\cup D$ is compact, and $\mathcal{A}$ is UGAS for $\mathcal{H}_2^{\text{ave}}$, then $\mathcal{A}$ is GPAS as $\varepsilon\to0^+$ for the HDS \eqref{eq:ms_orig_hybrid_sys} and the trajectories satisfy
\begin{align}\label{KLboundSGPAS223}
        |\theta(t,j)|_{\mathcal{A}}\leq \beta(|\theta(0,0)|_{\mathcal{A}},t+j)+ \nu, 
\end{align}
for all $(t,j)\in\text{dom}(\theta,\tau)$. \QEDB 
\end{corollary}
\begin{remark}\label{twotunableparameters}
In Theorem \ref{thm:avg_UAS_implies_org_PUAS}, it is assumed that $\mathcal{H}_2^{\text{ave}}$ renders UAS the set $\mathcal{A}$. However, in certain cases $\mathcal{H}_2^{\text{ave}}$ might depend on additional tunable parameters $\eta>0$, and $\mathcal{A}$ might only be SGPAS as $\eta\to0^+$, see \cite{PoTe17Auto,PovedaNaliAuto20,krilavsevic2023learning}. In this case, the SGPAS result of Theorem \ref{thm:avg_UAS_implies_org_PUAS} still applies as $(\varepsilon,\eta)\to0^+$, where the tunable parameters are now sequentially tuned, starting with $\eta$. This observation follows directly by modifying the proof of Theorem \ref{thm:avg_UAS_implies_org_PUAS} as in \cite[Thm. 7]{PovedaNaliAuto20}.  {An example  will be presented in Section \ref{subsectionswitching} in the context of slow switching systems}. \QEDB 
\end{remark}
%

\section{Applications to Hybrid Model-Free Control and Optimization}
\label{sec:switching_systems}
 {We present three different applications to illustrate the results of Section \ref{closenesstrajectories} in the context of model-free optimization and control.} In all our examples, we will consider a HDS of the form \eqref{eq:ms_orig_hybrid_sys}, where $f_{\varepsilon}$ is allowed to switch between a finite number of vector fields, some of which might not necessarily have a stabilizing second-order average HDS. Specifically, we consider systems where the main state $x$ evolves in an application-dependent closed set $M\subset\mathbb{R}^{n_1}$, under the following dynamics: 
\begin{align}
\dot{x}&= f_\varepsilon(x,z,\tau_1,\tau_2)= \sum_{k=1}^2\varepsilon^{k-2}\phi_{z_1,k}({x},\tau_1,\tau_2),\label{vaerepsilon}
\end{align}
where $z_1$ is a logic mode allowed to switch between values in the set $\mathcal{Q}=\{1,\dots,N\}$, where $N\in\mathbb{Z}_{\geq 2}$.  {The switching behavior can be either \emph{time-triggered} (modeled via the hybrid automaton \eqref{hybridautomaton}) or \emph{state-triggered} (modeled via a hysteresis-based mechanism incorporated into the sets $C$ and $D$). In both cases, Zeno behavior will be precluded by design.} Since for each $z_1\in \mathcal{Q}$ the functions $\phi_{z_1,k}$ will be designed to satisfy \thref{asmp:A1},  the \emph{2$^{\text{nd}}$-order average mapping} can be directly computed for each mode $z_1\in\mathcal{Q}$ using \eqref{fbareq}:
\begin{equation}\label{formulaaveragemap}
    \begin{aligned}
        \bar{f}_{z_1}&= \frac{1}{T_1T_2}\int_0^{T_1}\int_0^{T_2} \phi_{z_1,2}\,d\tau_2\,d\tau_1\\
        &+\frac{1}{T_1T_2}\int_0^{T_1}\int_0^{T_2}\frac{1}{2}\left[\int_0^{\tau_2} \phi_{z_1,1}\,ds_2,\phi_{z_1,1}\right]\,d\tau_2\,d\tau_1.
    \end{aligned}
\end{equation}
To guarantee that $\bar{f}_{z_1}$ is well-defined,  we let $\phi_{z_1,k}$ depend on a vector of frequencies $w=(w_1,w_2,\ldots,w_r)$, for some $r\in\mathbb{Z}_+$, which satisfies the following:
\vspace{0.1cm}
\begin{asmp}\label{assumptionfrequencies}
The frequencies satisfy $w_i\in\mathbb{Q}_{> 0}$ and $w_{i}\neq w_{j}$, $\forall~i\neq j\in\mathbb{N}$. \QEDB 
\end{asmp}

 {The conditions of Assumption \ref{assumptionfrequencies} are standard in vibrational control \cite{scheinker2017model} and extremum seeking \cite{KrsticBookESC,TanAndNesic2006Local,PoTe17Auto,suttner2022robustness,grushkovskaya2018class}}. 

\vspace{-0.1cm}
\subsection{Switching Model-Free Optimization and Stabilization}
\label{subsectionswitching}
In our first two applications, the set of logic modes satisfies $\mathcal{Q}=\mathcal{Q}_s\cup \mathcal{Q}_u$, where $\mathcal{Q}_s\cap \mathcal{Q}_u=\emptyset$, $\mathcal{Q}_s$ indicates ``stable'' modes, and $\mathcal{Q}_u$ indicates ``unstable'' modes. The logic mode $z_1:\text{dom}(z_1)\to\mathcal{Q}$ models a switching signal generated by the hybrid automaton \eqref{hybridautomaton} studied in Section \ref{secmotivationalexample}, with state $z=(z_1,z_2,z_3)$, which induces the ADT constraint \eqref{ADT} and the AAT constraint \eqref{ATT} on the switching, effectively ruling out Zeno behavior. Using the data $(C_z,\mathcal{T}_F,D_z,\mathcal{T}_G)$ of system \eqref{hybridautomaton}, the closed-loop switched system has the form of \eqref{eq:ms_orig_hybrid_sys}, with sets
\begin{equation}\label{flowjumpsetsexample}
C=M\times C_z,~~~D=M\times D_z,
\end{equation}
flow map $F_{\varepsilon}=\{f_{\varepsilon}\}\times \Phi$,  $f_{\varepsilon}$ given by \eqref{vaerepsilon}, $\Phi(z)=\mathcal{T}_F(z)$, and jump map:
\begin{equation}\label{overaljumpmapexample}
G(x,z)=\{x^+\}\times \mathcal{T}_G(z).
\end{equation}
 By construction, Assumption \ref{asmp:basic} is satisfied and the corresponding HDS $\mathcal{H}_2^{\text{ave}}$ can be directly obtained using \eqref{eq:ms_avg_hybrid_sys}:
\begin{subequations}\label{eq:switched_averaged_sys}
\begin{align}
(\bar{x},\bar{z})&\in C,\,\,\begin{cases}    
\left(\begin{array}{c}\dot{\bar{{x}}}\\\dot{\bar{z}}\end{array}\right)\in F(\bar{x},\bar{z}),
\end{cases} \\
(\bar{x},\bar{z})&\in D,\,\,\begin{cases}
\left(\begin{array}{c}\bar{{x}}^+\\{\bar{z}}^+\end{array}\right)\in G(\bar{x},\bar{z}),
\end{cases}
\end{align}
\end{subequations}
where $F(\bar{x},\bar{z})=\left\{
         \bar{f}_{z_1}(\bar{x})\right\}\times\mathcal{T}_{F}(\bar{z})$. 
\vspace{0.1cm}
We consider two different applications where this model is applicable. 

\vspace{0.1cm}
\subsubsection{Synchronization of Oscillators with  Switching Control Directions and Graphs} Consider a collection of $r\in\mathbb{Z}_{\geq2}$ oscillators evolving on the $r$-torus $M=\mathbb{S}^1\times\dots\times\mathbb{S}^1=\mathbb{T}^r\subset \mathbb{R}^{2r}$, which is an embedded (and compact) submanifold in the Euclidean space.  Let $x^i=(x^i_1,x^i_2)\in\mathbb{S}^1$ be the state of the $i^{\text{th}}$ oscillator, for $i\in\{1,\dots,r\}$, and denote its input by $u_i$.  {Moreover, let $\alpha:=(\alpha_1,\dots,\alpha_r)\in\mathcal{J}$, where $\mathcal{J}\subseteq\{+1,-1\}^r$, $|\mathcal{J}|=N_1\in\mathbb{N}$, such that $\alpha_i\in\{+1,-1\}$ for all $i$. The dynamics of the overall network of oscillators are given by
\begin{align}\label{eq:synchronization_exmp_org_sys_full_state}
\begin{pmatrix} \dot{x}^1 \\ \vdots \\ \dot{x}^r\end{pmatrix} = \begin{pmatrix} (1+\alpha_1 u_1)\,S x^1 \\ \vdots \\ (1+\alpha_r u_r)\,S x^r\end{pmatrix},~~S:= \begin{pmatrix} 0 & 1 \\-1 & 0\end{pmatrix}.
\end{align}
where $\alpha$ is the unknown vector of control directions.} Let $x=(x^1,\dots,x^r)\in M$, and note that $M$ is invariant under the dynamics in \eqref{eq:synchronization_exmp_org_sys_full_state}. The connectivity between oscillators is dictated by the network topology which, at any given instant, is defined by a unique element from a collection of $N_2\in\mathbb{Z}_{\geq1}$ graphs $\{\mathcal{G}_1,\dots,\mathcal{G}_{N_2}\}$, $\mathcal{G}_k=(\mathcal{V},\mathcal{E}_k)$, where $\mathcal{V}=\{1,\dots,r\}$ is the set of nodes representing the oscillators, and $\mathcal{E}_k$ is the edge set representing the interconnections. To simplify our presentation, we make the following assumption.

\vspace{0.1cm}
\begin{asmp}\thlabel{asmp:ex_synchronization_A1}
For each $k\in\{1,2,\ldots,N_2\}$ the graph $\mathcal{G}_k$ is connected and undirected.  \QEDB 
\end{asmp}

\begin{figure}[t]
    \centering
    \includegraphics[width=0.475\textwidth]{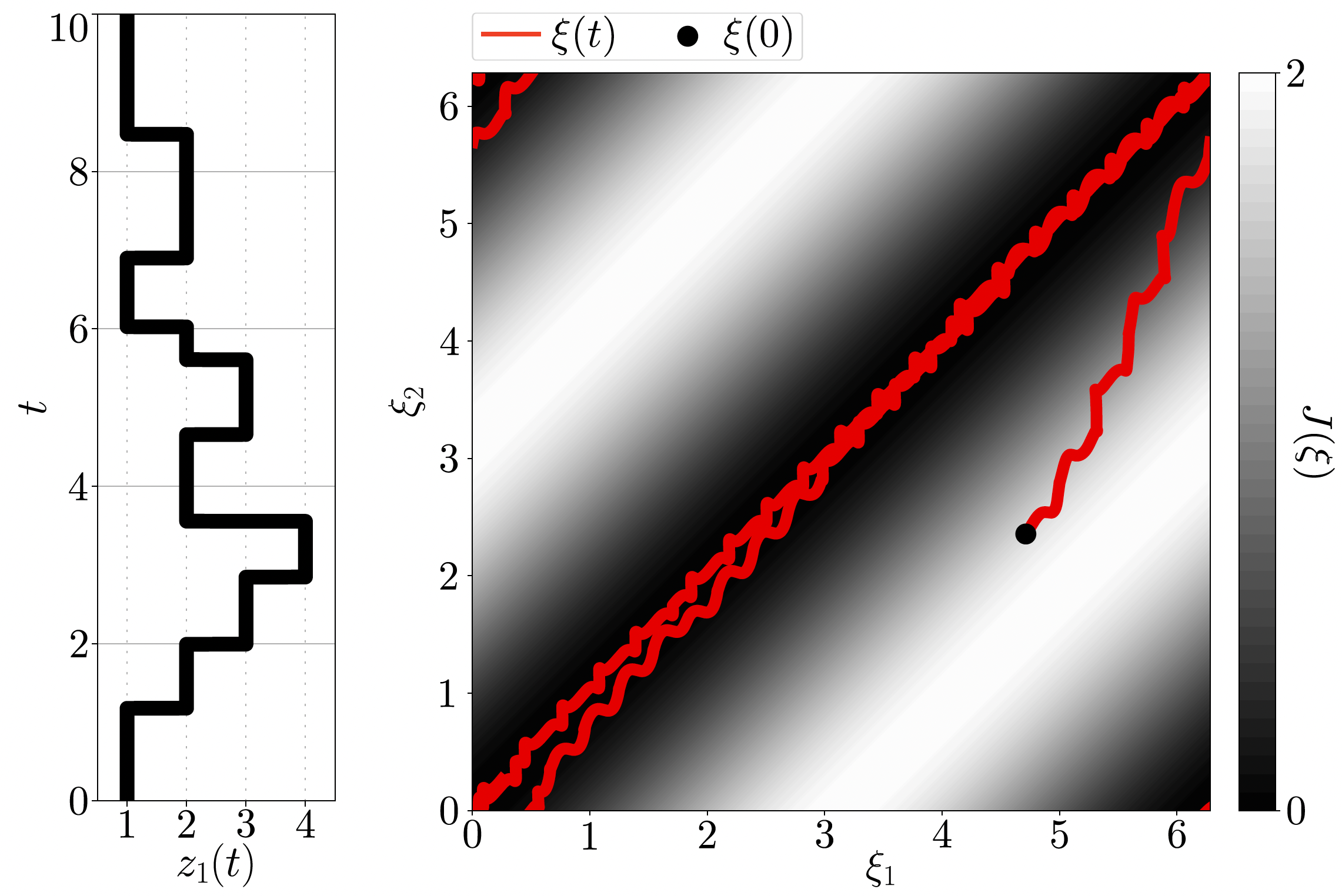}
    \caption{\small{ {Simulation results of the first scenario in Example 1. The figure on the right is a flat embedding of the torus $\mathbb{T}^2\subset\mathbb{R}^4$ into the plane $\mathbb{R}^2$. Synchronization is achieved on the submanifold $\mathcal{S}$ corresponding to the diagonal line.}}}
\label{fig:synchronization_with_switching_topologies_exmp_1}
    \vspace{-0.3cm}
\end{figure}

\begin{figure}[t]
    \centering
    \includegraphics[width=0.475\textwidth]{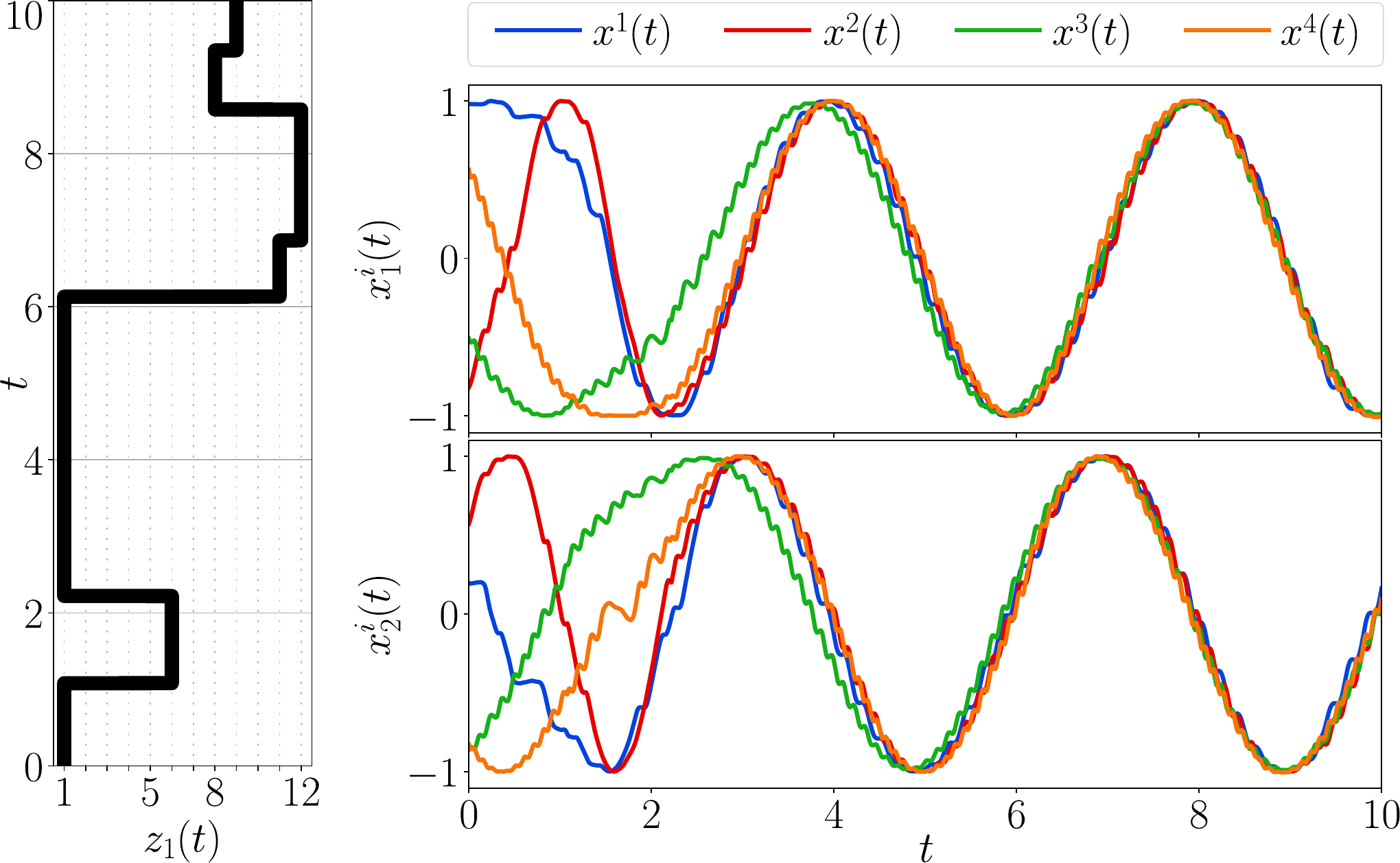}
    \caption{\small{ {{Simulation results of the second scenario considered in Example 1.}}}}
\label{fig:synchronization_with_switching_topologies_exmp_2}
\vspace{-0.3cm}
\end{figure}

 {Our goal is to synchronize the oscillators under switching communication topologies $\mathcal{G}_k$ and unknown control directions $\alpha$.} To solve this problem,  {we consider the feedback law
\begin{align}\label{eq:synchronization_exmp_feedback_law}
    u_{i,k}(x,\tau_2)&= \varepsilon^{-1}\sqrt{2w_i}\cos\left(w_i \tau_2 + \kappa J_{i,k}(x)\right)\kappa^{-\frac{1}{2}},
\end{align}    
}where $J_{i,k}(x)= \frac{1}{2}\sum_{j\in\mathcal{N}^i_{k}}\lvert x^i-x^j\lvert^2$, and $\kappa\in\mathbb{R}_{>0}$ is a tuning parameter. In \eqref{eq:synchronization_exmp_feedback_law}, the set $\mathcal{N}_{k}^i=\{j\in\mathcal{V}: (i,j)\in\mathcal{E}_k\}\subset\mathcal{V}$ indicates the neighbors of the $i^{\text{th}}$ agent according to the topology of the graph $\mathcal{G}_k$.  We consider a model where both the network topology $\mathcal{G}_k$ and the control directions $\alpha$ are allowed to switch in \eqref{eq:synchronization_exmp_org_sys_full_state} and \eqref{eq:synchronization_exmp_feedback_law}. To capture this behavior, let $\mathcal{Q}=\{1,\dots,N_1N_2\}$, and fix a choice of all possible bijections $\varrho:\mathcal{Q}\ni z_1 \mapsto(k(z_1),\alpha(z_1))\in\{1,\dots,N_2\}\times\mathcal{J}$. In this way, for each $z_1\in \mathcal{Q}$ there is a particular  control direction $\alpha(z_1)$ and graph $\mathcal{G}_{k(z_1)}$ acting in \eqref{eq:synchronization_exmp_org_sys_full_state} and \eqref{eq:synchronization_exmp_feedback_law}.
To study the overall dynamics of the system, we define the coordinates
\begin{align*}
    x^i_1&= -\cos(\xi_i), & x^i_2&= \sin(\xi_i), & \forall i&\in\mathcal{V},
\end{align*}
where $\xi_i\in[0,2\pi)$, and $\xi:=(\xi_1,\dots,\xi_r)$, which lead to the following polar representation of  \eqref{eq:synchronization_exmp_org_sys_full_state}:
\begin{align}\label{eq:synchronization_exmp_org_sys_full_state_polar}
\begin{pmatrix} \dot{\xi}_1 \\ \vdots \\ \dot{\xi}_r\end{pmatrix} = \begin{pmatrix} 1+\alpha_1 u_{1,k}(\xi,\tau_2) \\ \vdots \\ 1+\alpha_r u_{r,k}(\xi,\tau_2)\end{pmatrix},
\end{align}
\begin{figure}
    \centering
    \includegraphics[width=0.45\textwidth]{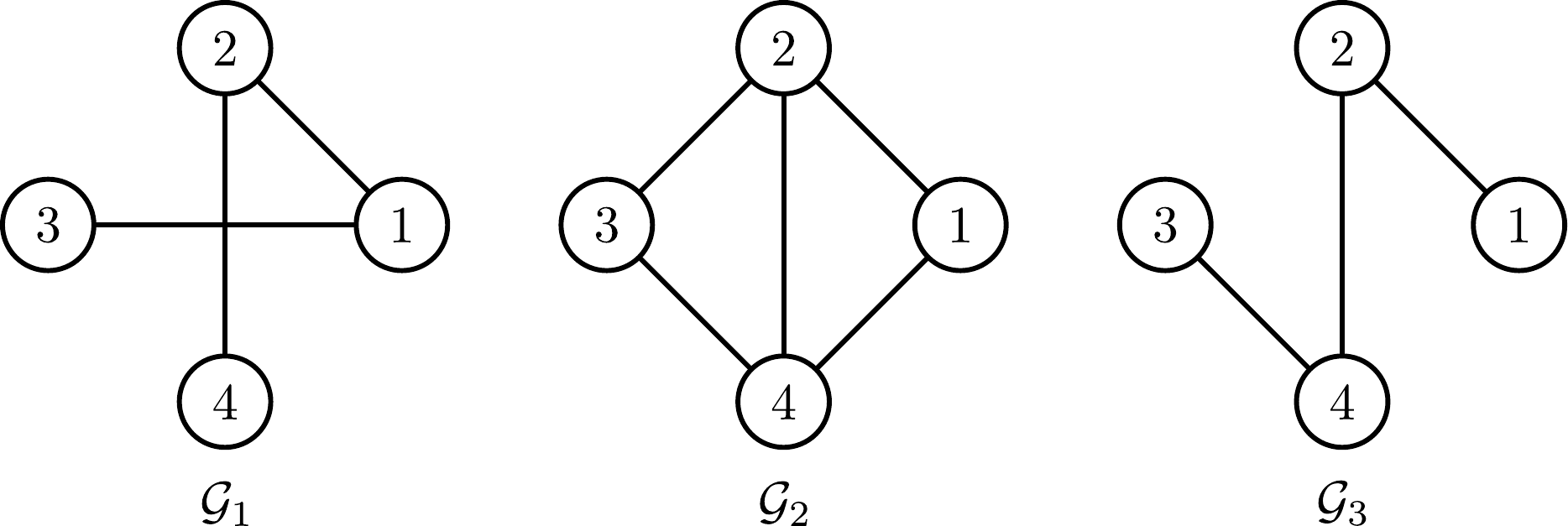}
    \caption{\small{Network topologies in the second case of Example 1.}}
    \label{fig:network_topologies}
    \vspace{-0.3cm}
\end{figure}
where the feedback law $u_{i,k}$ becomes:
\begin{align}\label{eq:synchronization_exmp_feedback_law_polar}
    u_{i,k}(\xi,\tau_2)&= \varepsilon^{-1}\sqrt{2w_i}\cos(w_i\tau_2+\kappa J_{i,k}(\xi))\kappa^{-\frac{1}{2}}, 
\end{align}
with distance function $J_{i,k}(\xi) = \sum_{j\in\mathcal{N}^i_k}(1-\cos(\xi_i-\xi_j))$. Using \eqref{formulaaveragemap}, the $2^{\text{nd}}$-order average mapping is given by
\begin{align}\label{eq:synchronization_exmp_avg_sys_full_state_polar}
\dot{\bar{\xi}} = \mathbbm{1} - \nabla V_{z_1}(\bar{\xi}),
\end{align} 
where $\mathbbm{1}=(1,1,\ldots,1)$, $V_{z_1}(\bar{\xi})= \sum_{i=1}^{r} J_{i,k(z_1)}(\bar{\xi})$, and
%
%
%
\begin{align*}
    \nabla V_{z_1}(\bar{\xi})&= \begin{pmatrix}
        \sum_{j\in\mathcal{N}^1_{k(z_1)}}\sin(\bar{\xi}_1-\bar{\xi}_j)\\
        \vdots\\
        \sum_{j\in\mathcal{N}^r_{k(z_1)}}\sin(\bar{\xi}_r-\bar{\xi}_j)
    \end{pmatrix},
\end{align*}
which is the well-known Kuramoto model over the graph $\mathcal{G}_{k(z_1)}$ \cite{durr2013examples,DurrManifold}. The subset $\mathcal{S}$ that characterizes synchronization has the polar coordinate representation $   \mathcal{S}= \{\xi\in [0,2\pi)^r: \,\xi_1 = \cdots = \xi_r\}.$ Now observe that the vector field in (\ref{eq:synchronization_exmp_org_sys_full_state}) under the feedback law (\ref{eq:synchronization_exmp_feedback_law}), or equivalently the vector field in (\ref{eq:synchronization_exmp_org_sys_full_state_polar}) under the feedback law (\ref{eq:synchronization_exmp_feedback_law_polar}), has the same structure as the vector fields $f_\varepsilon$ in \eqref{vaerepsilon} with trivial dependence on the intermediate timescale $\tau_1=\varepsilon^{-1}t$. Moreover, using  \eqref{hybridautomaton} with $\mathcal{Q}_u=\emptyset$ to generate the switching signal $z_1$, the closed-loop system has the form \eqref{eq:ms_orig_hybrid_sys} with $C$ and $D$ given by \eqref{flowjumpsetsexample}, and $G$ given by \eqref{overaljumpmapexample}. For this system, we study stability of the set $\mathcal{A}=\mathcal{S}\times C_z$ using two tunable parameters (c.f. Remark \ref{twotunableparameters}).   {The following proposition establishes a novel ``model-free'' synchronization result for oscillators under switching graphs and control directions via averaging.}
\begin{prop}\label{prop:switching_topology_es}
Under Assumptions \ref{assumptionfrequencies}-\ref{asmp:ex_synchronization_A1}, the HDS \eqref{eq:ms_orig_hybrid_sys} renders the set $\mathcal{A}$ SGPAS as $(\varepsilon,\eta_1)\rightarrow 0^+$ with some basin of attraction $\mathcal{B}_{\mathcal{A}}$. \QEDB 
\end{prop}
%
%
%
%
\begin{example}
To illustrate the above result, we consider two scenarios. In the first scenario, $\mathcal{G}_k$ is static and the control directions are switching. In this case, we consider two oscillators (i.e., $r=2$), $N_2=1$ and $J_{1,2}(\xi)=J_{2,1}(\xi)=J(\xi)=1-\cos(\xi_1-\xi_2)$. We let $N_1=4$, $\mathcal{J}=\{(+1,+1),(-1,+1),(+1,-1),(-1,-1)\}$, and we consider the bijection $\varrho$ that satisfies $\varrho(1)= (1,(+1,+1))$, $\varrho(2)= (1,(-1,+1))$, $\varrho(3)= (1,(+1,-1))$, $\varrho(4)= (1,(-1,-1))$. For the numerical simulation results, shown in Fig. \ref{fig:synchronization_with_switching_topologies_exmp_1}, we used  {$\varepsilon=1/\sqrt{10\pi}$, $\kappa = 10$, $w_1=1$, $w_2=2$, and a switching mode $z_1$ generated by \eqref{hybridautomaton} with $\eta_1=2.5$, $N_{\circ}=1$, and arbitrary $\eta_2$ and $T_{\circ}$.} As observed in the figure, (local) practical synchronization is achieved. 

In the second scenario, we incorporate switching graphs $\mathcal{G}_k$. We let $r=4$, $N_1=4$, $N_2=3$, and we consider the set $\mathcal{J}= \{(+1,+1,-1,+1), (-1,+1,+1,+1),(-1,+1,-1,-1),\\(-1,-1,+1,+1)\}$, and the collection of graphs $\{\mathcal{G}_1,\mathcal{G}_2,\mathcal{G}_3\}$ shown in Fig. \ref{fig:network_topologies}. 
For the numerical simulation results, shown in Fig. \ref{fig:synchronization_with_switching_topologies_exmp_2}, we used  {$\varepsilon=1/\sqrt{10\pi}$, $\kappa = 10$, $w_1=1$, $w_2=4/3$, $w_3=5/3$, $w_4=2$, and a switching mode $z_1$ generated by \eqref{hybridautomaton} with $\eta_1=1.5$, $N_{\circ}=1$, and arbitrary $\eta_2$ and $T_{\circ}$.}
As shown in Fig. \ref{fig:synchronization_with_switching_topologies_exmp_2}, (local) practical synchronization is also achieved despite simultaneous switches happening in the network topology and the control directions. \QEDB  
\end{example}

\vspace{0.1cm}
\subsubsection{Lie-Bracket Extremum Seeking under Intermittence and Spoofing}
\label{intermittentliebracket}
 {Consider the following control-affine system  
\begin{align}\label{eq:intermittent_es_ex_org_sys}
    \dot{x}= \sum_{i=1}^r b_i(x,\tau_1) u_{i,z_1}(x,\tau_2),~~~r\in\mathbb{Z}_{\geq1},
\end{align}
where the goal is to steer the state $x\in\mathbb{R}^n$ towards the set of solutions of the optimization problem
\begin{equation}\label{ESproblemopti}
\tcb{\min_{x\in M}}~~J(x),
\end{equation}
where $M\subset\mathbb{R}^n$ is the parameter search space, and $J$ is a continuously differentiable cost function whose mathematical form is unknown, but which is available via measurements or evaluations. Problem \eqref{ESproblemopti} describes a standard extremum seeking (ES) problem \cite{DurrLieBracket,KrsticBookESC,PoTe17Auto,scheinker2012minimum}, where exact knowledge of the functions $b_i$ is in general not required.}  {However, unlike traditional ES, we aim to achieve convergence to the solutions of \eqref{ESproblemopti} despite sporadic failures in accessing measurements of $J$ and intermittent cost measurements affected by \emph{malicious external spoofing} designed to destabilize the system, as in the motivational example of Section \ref{secmotivationalexample}, see also \cite{labar2022extremum,galarza2021extremum}.} 

To study this scenario, we consider the following mode-dependent Lie-bracket-based control law
\begin{align}\label{eq:intermittent_es_ex_feedback_law}
    u_{i,{z_1}}(x,\tau_2) &= \varepsilon^{-1}\sqrt{2w_i}\cos\left(w_i \tau_2 + \kappa (z_1-2) J(x)\right)\kappa^{-\frac{1}{2}},
\end{align}

\vspace{-0.3cm}
\noindent 
where $\kappa\in\mathbb{R}_{>0}$ is a tuning parameter and $z_1\in\mathcal{Q}=\{1,2,3\}$ is the logic mode generated by the hybrid automaton \eqref{hybridautomaton} . Note that  \eqref{nonholonomicvehiclea1} is a particular case of \eqref{eq:intermittent_es_ex_org_sys}.  Using equation \eqref{formulaaveragemap} to study system \eqref{eq:intermittent_es_ex_org_sys} under the feedback law \eqref{eq:intermittent_es_ex_feedback_law}, we obtain:
\begin{align*}
    \bar{f}_{z_1}(\bar{x})&= (2-z_1)P(\bar{x})\,\nabla J(\bar{x}),
\end{align*}
where the matrix-valued mapping $P$ is given by
\begin{align*}
    P(\bar{x}) =\frac{1}{T_1} \sum_{i=1}^r\int_0^{T_1} b_i(\bar{x},\tau_1)\otimes b_i(\bar{x},\tau_1)\,d\tau_1.
\end{align*}
We make the following regularity assumptions on $J$ and  $P$:
\vspace{0.1cm}
\begin{asmp}\thlabel{asmp:ex_2_A1}
The following holds:
\begin{enumerate}
\item There exists $\lambda_{P}>0$ and $M_{P}>0$ such that $|P(x)x|\leq M_P|x|, \, x^\intercal P(x) x \geq \lambda_{P}|x|^2, \, \forall x\in\mathbb{R}^n.$
\item  The cost function $J$ is $\mathcal{C}^1$, strongly convex with strong-convexity parameter $\mu>0$, and minimizer $x^{\star}$, and \tcb{there exists $L_J>0$} such that $|\nabla J(x)-\nabla J(x')|\leq L_{J}|x-x'|$, for all $x,x'\in\mathbb{R}^n$. \QEDB 
\end{enumerate}
\end{asmp}

\vspace{0.1cm}
\begin{remark}
\tcb{Item 1) in Assumption \ref{asmp:ex_2_A1} can be seen as a cooperative persistence of excitation condition on the vector fields $b_i$, see \cite{abdelgalil2022recursive} for similar conditions. Item 2) is common in ES \cite{KrsticBookESC}, and covers quadratic functions}. \QEDB 
\end{remark}

\vspace{0.1cm}
 {The following result formalizes the discussions of Section \ref{secmotivationalexample}, and establishes a novel resilience result for Lie-bracket ES algorithms under spoofing and intermittent feedback}: 

\vspace{0.1cm}
\begin{prop}\thlabel{prop:ex2_es_intermittent}
Under Assumptions \ref{assumptionfrequencies} and  \ref{asmp:ex_2_A1}, for $\eta_1>0$ there exist $\eta_2^*>0$ such that for all $\eta_2\in(0,\eta_2^*)$ the HDS \eqref{eq:ms_orig_hybrid_sys} renders the set $\mathcal{A}=\{x^{\star}\}\times C_z$ SGPAS as $\varepsilon\to0^+$. \hfill $\QEDB$
\end{prop}
\subsection{Switched Global Extremum Seeking on Smooth Compact Manifolds}
We now leverage the tools presented in this paper to solve  {\emph{global} ES problems of the form \eqref{ESproblemopti} on sets $M$ that describe a class of smooth, boundaryless, compact manifolds.} This application illustrates the use of Corollary \ref{corollary1}. As thoroughly discussed in the literature \cite{mayhew2011topological, TrackingManifold, ochoa2025robust,bernuau2013retraction}, for such problems, no smooth algorithm with global stability bounds of the form \eqref{KLboundSGPAS223} can be employed as the average system due to inherent topological obstructions. Here, we show how to remove this obstruction by using adaptive switching between multiple ES controllers in problems that satisfy the following assumption:

\vspace{0.1cm}
\begin{asmp}\thlabel{asmp:mfd_ex_A1}
The following holds:
\begin{enumerate}[(a)]
\item The set $M$ is a smooth, embedded, connected, and compact submanifold without boundary, endowed with a Riemannian structure by the metric $\langle\cdot,\cdot\rangle_M:T_xM\times T_xM\rightarrow \mathbb{R}$ induced by the metric $\langle\cdot,\cdot\rangle:\mathbb{R}^n\times \mathbb{R}^n\rightarrow \mathbb{R}$ of the ambient Euclidean space $\mathbb{R}^n$.
\item The function $J$ is smooth on an open set $U\subset\mathbb{R}^n$ containing $M$, has a unique minimizer $x^{\star}\in M$ satisfying $J(x^{\star}) < J(x)$, $\forall x\neq x^\star\in M$, and there exists a known number $ \bar{J}\in(0,\infty)$ such that for all $x\in \mathcal{L}_{\bar{J}} = \{x\in M:\, J(x)< J(x^{\star}) + \bar{J}\}$ we have $\nabla_M J(x) =0 \iff x = x^{\star}$, where $\nabla_M J$ is the orthonormal projection (defined by the metric $\langle\cdot,\cdot\rangle$) of the gradient $\nabla J$ onto the tangent space $T_x M$, $\forall x\in M$. \QEDB 
\end{enumerate}
\end{asmp}

\vspace{0.1cm}
 {While \thref{asmp:mfd_ex_A1} suffices to guarantee the local, or at best, almost-global uniform practical asymptotic stability of the minimizer $x^{\star}$ for a smooth Lie Bracket-based extremum seeking controller on $M$ \cite{durr2013examples,DurrManifold}, its (smooth) average system will necessarily have more than one critical point in $M$ due to the topology of the manifold \cite{mayhew2011topological,Mayhew10Thesis}. To remove this obstruction to achieve uniform global ES, we consider a hybrid algorithm that implements state-triggered switching between certain functions introduced in the following definition:}

\vspace{0.1cm}
\begin{definition}\label{defmanifolds}
Let $N\in\mathbb{N}_{\geq 1}$, and $\mathcal{Q}=\{1,\dots,N\}$. A family of functions $\{\tilde{J}_q\}_{q\in\mathcal{Q}}$ is said to be a $\delta$-gap synergistic family of functions subordinate to $J$ on $M$ if:
\begin{enumerate}
    \item  $\forall q\in\mathcal{Q}$, $\tilde{J}_q: U \rightarrow \mathbb{R}$ is smooth, its unique minimizer on $M$ is $x^{\star}$, and has a finite number of critical points  in $M$, i.e. $|\text{Crit }\tilde{J}_q|\in\mathbb{N}_{\geq 1}$, where $\text{Crit }\tilde{J}_q:=\{x\in M:\,\nabla_M \tilde{J}_q(x)=0\}$.
    \item $\forall q\in\mathcal{Q}$, $J(x) < \bar{J}+J(x^{\star}) \implies \tilde{J}_q(x) = J(x)$, i.e. the family $\{\tilde{J}_q\}$ agrees with $J$ on the neighborhood $\mathcal{L}_{\bar{J}}$ of the minimizer $x^{\star}$,
    \item $\exists~\delta\in(0,\infty)$ such that
    \begin{align*}
        \delta <\Delta^*:=\min_{\substack{q_1\in\mathcal{Q}\\x\in\text{Crit }\tilde{J}_{q_1}\backslash\{x^{\star}\}}}\max_{q_2\in\mathcal{Q}}\,(\tilde{J}_{q_1}(x)-\tilde{J}_{q_2}(x)). ~~~\vspace{0.1cm}\QEDB
    \end{align*} 
\end{enumerate}
\end{definition}
\begin{figure}[t]
        \centering    \includegraphics[width=0.46\textwidth]{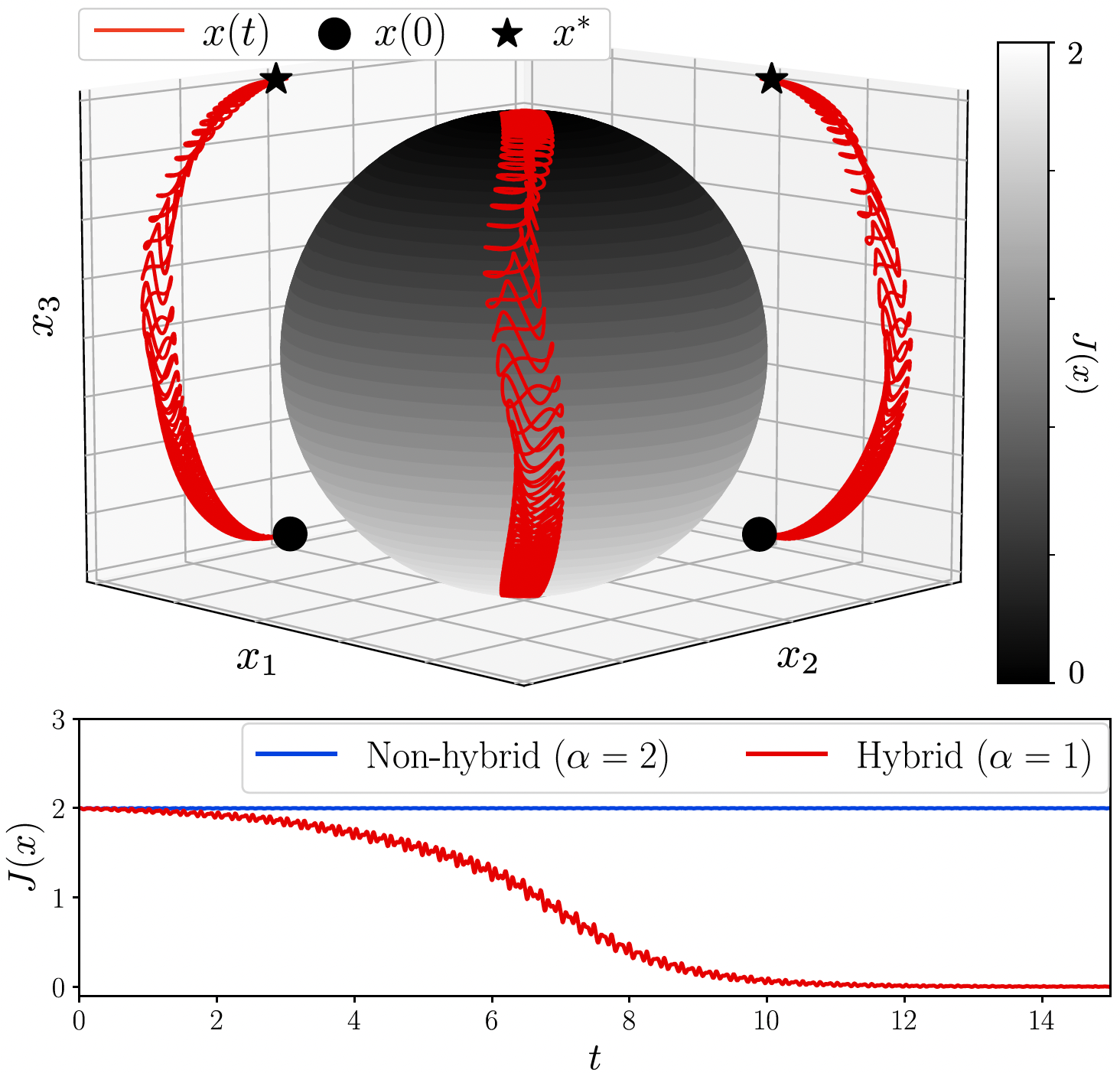}
        \caption{\small{ {{Simulation results of the first scenario in \thref{exmp:optimization_on_manifolds_example}: (top) the trajectory on $\mathbb{S}^2$ for $\alpha = 1$ along with two projections,  {and (bottom) the evolution of the cost along the trajectories.}}}}}\label{fig:optimization_on_manifolds_example_trajectory_visualization}
        \vspace{-0.2cm}
\end{figure}
In words, the constructions in Definition \ref{defmanifolds} involve families of diffeomorphisms that ``warp the manifold $M$ in sufficiently many ways so as to distinguish the minimizer $x^{\star}$ from other critical points of the function $J$". The reader is referred to \cite{Strizic:17_CDC} and \cite{ochoa2025robust} for the full details (and examples) of these constructions, which are leveraged in \cite{ochoa2025robust} for the design of algorithms using classic averaging.  {We stress that constructing these diffeomorphisms does not necessarily require exact knowledge of the mathematical form of $J$ (although knowledge of $M$ is required), but only a qualitative description of $J$ that allows to identify the threshold $\bar{J}$ below which the only critical point corresponds to $x^{\star}$. Such knowledge is reasonable for most practical applications where the optimal value $J(x^{\star})$ is known to lie within a certain bounded region, or be equal to certain value.}  {A similar assumption was considered in \cite{suttner2017exponential}. Under such qualitative knowledge, our goal is to attain \emph{global} ES from all initial conditions, ensuring that \eqref{ESproblemopti} is solved regardless of where the algorithm's trajectories are initialized or land after a sudden disturbance.}

To design hybrid Lie-bracket-based algorithms with global stability properties, we consider dynamics of the form \eqref{eq:intermittent_es_ex_org_sys}, where, for simplicity we allow the functions $b_i$ to be independent of $\tau_1$. We formalize this in the following assumption:

\begin{asmp}\thlabel{asmp:mfd_ex_A2}
The following holds for problem \eqref{ESproblemopti}:
\begin{enumerate}[(a)]
\item There exists a family of smooth vector fields $\{b_1,\dots,b_r\}$ on $M$, for which the operator $P:M\times \mathbb{R}^n\ni(x,v)\mapsto P(x)[v]\in T_x M$ defined by $$P(x) [v]= \sum_{i=1}^r\langle b_i(x),v\rangle\,b_i(x),$$ is such that $P(x) [v] = 0 \iff v\in \left(T_xM\right)^\perp.$
\item  There exists a $\delta$-gap synergistic family of functions $\{\tilde{J}_q\}_{q\in\mathcal{Q}}$ subordinate to $J$ on $M$, which are available for measurement.  \QEDB 
   
\end{enumerate}   
\end{asmp}
\vspace{0.1cm}

\tcb{Assumption \ref{asmp:mfd_ex_A2} is natural in the setting we consider herein. Indeed, item (a) in Assumption \ref{asmp:mfd_ex_A2} is a standard assumption in the context of Lie-bracket extremum seeking on manifolds, see e.g. \cite[Assumption 1]{DurrManifold}. On the other hand, item (b) is a standard assumption in the context of \emph{robust} global stabilization problems on manifolds via hybrid control, see e.g. \cite{ochoa2025robust,mayhew2011topological,casau2024robust}.}

To solve \eqref{ESproblemopti}, we consider the HDS \eqref{eq:ms_orig_hybrid_sys}, with $f_\varepsilon$ given by \eqref{eq:intermittent_es_ex_org_sys}, and
%
%
 {the following hybrid feedback law
\begin{align*}
    u_{i,z}(x,\tau_2):=\varepsilon^{-1}\sqrt{2w_i}\cos(\kappa\,\tilde{J}_z(x)+w_i\tau_2)\kappa^{-\frac{1}{2}},
\end{align*}}
\!\!\!where $\kappa\in\mathbb{R}_{>0}$ is a tuning gain, $z\in\mathcal{Q}$ is a logic state, the sets $C$ and $D$ are given by
\begin{align*}
    C&:=\big\{(x,z)\in M\times\mathcal{Q}:\, \tilde{J}_z(x)-\min_{\tilde{z}\in\mathcal{Q}}\tilde{J}_{\tilde{z}}(x) \leq \delta\big\},\\
    D&:=\big\{(x,z)\in M\times\mathcal{Q}:\, \tilde{J}_z(x)-\min_{\tilde{z}\in\mathcal{Q}}\tilde{J}_{\tilde{z}}(x) \geq \delta\big\},
\end{align*}
and the flow map in \eqref{flowmapstructure} uses $\Phi(z) = \{0\}$. Finally, the jump map in \eqref{secondcaseHDS} is defined as 
\begin{align*}
        G(x,z)= \{x\}\times
        \big\{z\in\mathcal{Q}:~\tilde{J}_{z}(x)=\min_{\tilde{z}\in\mathcal{Q}}\tilde{J}_{\tilde{z}}(x)\big\}.
\end{align*}
Using equation \eqref{formulaaveragemap} with $z=z_1$, we obtain:
\begin{equation}\label{averagemappingmanifold}
    \bar{f}_z(\bar{x})= -P(\bar{x})[\nabla \tilde{J}_{z}(\bar{x})].
\end{equation}
 {We can now state the following \emph{global} result for Lie-bracket ES systems on smooth compact manifolds:}
\vspace{0.1cm}
\begin{prop}\label{prop:global_es_example}
    Under Assumptions \ref{assumptionfrequencies}, \ref{asmp:mfd_ex_A1}, and \ref{asmp:mfd_ex_A2}, the set $\mathcal{A}=\{x^*\}\times\mathcal{Q}$ is GPAS as $\varepsilon\to0^+$ for the HDS  \eqref{eq:ms_orig_hybrid_sys}. \QEDB 
\end{prop}

\vspace{0.1cm}
\begin{figure}[t]
        \centering    \includegraphics[width=0.45\textwidth]{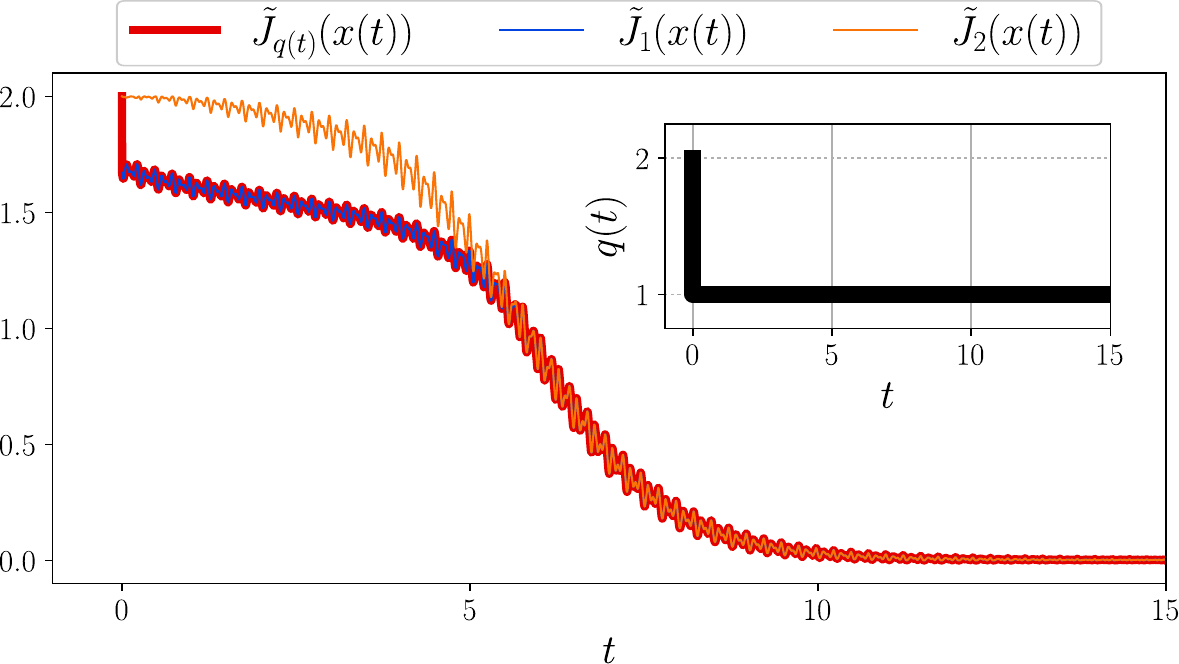}
        \caption{\small{ {{Simulation results of the second scenario in \thref{exmp:optimization_on_manifolds_example}.}}}}\label{fig:optimization_on_manifolds_example_numerical_simulation}
        \vspace{-0.4cm}
\end{figure}

\begin{example}\thlabel{exmp:optimization_on_manifolds_example}
    To illustrate the result of Proposition \ref{prop:global_es_example}, let $M=\mathbb{S}^2=\{x\in\mathbb{R}^3: |x|^2=1\}$ be  {the unit sphere}, and consider the cost $J(x)=1 - \langle x,e_3\rangle$. We define the vector fields $b_i$ for $i\in\{1,2,3\}$ by $b_i(x)= e_i\times x$, where $e_i$ is the unit vector in $\mathbb{R}^3$ with zero elements except the $i^{\text{th}}$-element which is set to $1$. The function $J$ has two critical points where $\nabla_{\mathbb{S}^2} J$ vanishes: a critical point $x^{\star}=(0,0,1)$ that corresponds to the minimum value $J(x^{\star})=0$, and a critical point $x^\sharp=(0,0,-1)$ that corresponds to the maximum value $J(x^\sharp)=2$. Let $\mathcal{Q}=\{1,2\}$, and define the family of functions $\{\tilde{J}_1,\tilde{J}_2\}$ by $\tilde{J}_q(x)= J\left(\Phi_q(x)\right)$, where the maps $\Phi_q$ are defined as follows:
    \begin{align*}
        \Phi_q(x)&:=\begin{cases}
            x, & J(x) \leq \alpha\\
            \exp\left(\frac{k_q}{\sqrt{2}}(J(x)-1)^2 (\hat{e}_1+\hat{e}_2)\right)x, & J(x)> \alpha
        \end{cases},
    \end{align*}
    and where $k_1=+\frac{1}{2}, \,k_2=-\frac{1}{2}$, $\alpha\in\{1,2\}$, $\exp:\mathbb{R}^{3\times 3}\rightarrow GL(3,\mathbb{R})$ is the matrix exponential where $GL(3,\mathbb{R})$ is the group of real invertible $3\times 3$-matrices, and $\hat{e}_i$ is the skew-symmetric matrix associated with the unit vector $e_i$. It can be verified, see \cite{ochoa2025robust}, that, when $\alpha = 1$, the family $\{\tilde{J}_1,\tilde{J}_2\}$ is a $\delta$-gap synergistic family of functions subordinate to $J$ on $\mathbb{S}^2$ for any $\delta\in(0,\frac{1}{4})$. We use $\alpha=2$ to emulate the situation in which no switching takes place since, when $\alpha = 2$, $J(x)=\tilde{J}_1(x)=\tilde{J}_2(x)$, for all $x\in\mathbb{S}^2$.  {To generate the numerical results, we used $\delta=\frac{1}{5}$, $w_1=2,\,w_2=3,\,w_3=1$, $\kappa=4$, and $\varepsilon=1/\sqrt{8\pi}$, which completely defines the HDS \eqref{eq:ms_orig_hybrid_sys}. We simulate two scenarios. In the first scenario, a small \emph{adversarial input}, tailored to (locally) stabilize the critical point $x^\sharp$ for an ES algorithm \textit{without switching}, is added to the vector field $f_\varepsilon$, i.e., we use $\tilde{f}_\varepsilon = f_\varepsilon + f_e$, where $|f_e|<0.1$. For details on the construction of $f_e$, we refer the reader to the appendix in the extended manuscript \cite{abdelgalil2023lie}. The simulation results of the first scenario for $\alpha=1$ and $\alpha=2$, starting from the initial condition $x(0,0)\approx(-0.11,0.11,-0.98)$, which is nearby $x^\sharp$, are shown in Fig. \ref{fig:optimization_on_manifolds_example_trajectory_visualization}. As observed, without switching ($\alpha = 2$), the small disturbance $f_e$ effectively traps the non-hybrid ES algorithm in the vicinity of the problematic critical point $x^\sharp$, despite the fact that $x(0,0)\neq x^\sharp$ and that the system is persistently perturbed by the dithers. On the other hand, as predicted by Proposition \ref{prop:global_es_example}, the hybrid ES algorithm ($\alpha =1$) renders the set $\mathcal{A}$ GPAS as $\varepsilon\rightarrow 0^+$}. 
    
     {In the second scenario, a similarly constructed adversarial input, tailored to (locally) stabilize the problematic critical points of the functions $\tilde{J}_1$ and $\tilde{J}_2$, \emph{in the absence of switching}, is added to $f_\varepsilon$. The simulation results for $\alpha=1$, with $q(0,0)=2$ and starting from $x(0,0)\approx(-0.286, 0.286, -0.914)$, which coincides with the problematic critical point of $\tilde{J}_2$, are shown in Fig. \ref{fig:optimization_on_manifolds_example_numerical_simulation}. It can be observed that a state jump from $q=2$ to $q^+=1$ is trigger to allow the ES algorithm escape the problematic critical point of $\tilde{J}_2$. The detailed construction of the adversarial signals can be found in the Appendix of the Supplemental Material.}  {All computer codes used to generate the simulations of this paper are available at \cite{HOA-HDS-Github}.}  \QEDB 
    %
\end{example}
%
%
    
\section{Analysis and Proofs}
\label{sec:proofs}
 {In this section, we present the proofs of our main results.}

\vspace{-0.3cm}
\subsection{Proof of Theorem \ref{thm:closeness_of_trajs}}  We follow similar ideas as in standard averaging theorems \cite[Thm. 10.4]{khalil}, \cite[Thm.1]{TeelNesicAveraging}, but using a recursive variable transformation applied to the analysis of the HDS.   

Fix $K_0\subset\mathbb{R}^n$, $T>0$ and $\rho>0$, and let Assumption \ref{asmp:A1} generate $\delta^{\star}>0$. For any $\delta\in(0,\delta^{\star}]$, consider the $\delta$-inflation of $\mathcal{H}_{2}^{\text{ave}}$, given by
\begin{equation}\label{eq:ms_inflated_avg_hybrid_sys}
\mathcal{H}_{2,\delta}^{\text{ave}}:~~~\left\{\begin{array}{ll}
\bar{\theta}\in C_{\delta}, &   
~~\dot{\bar{\theta}}\in F_{\delta}\left(\bar{\theta}\right)
\vspace{0.2cm}\\
\bar{\theta}\in D_{\delta}, &     
\bar{\theta}^+\in G_{\delta}\left(\bar{\theta}\right)
\end{array}\right.,
\end{equation}
where the data $(C_{\delta},F_{\delta},D_{\delta},G_{\delta})$ is constructed as in Definition \ref{def:inflatedHDS} from the data of $\mathcal{H}_{2}^{\text{ave}}$ in \eqref{defaveragehybridsystem}. Fix $\delta\in(0,\delta^{\star}]$ such that for any solution $\tilde{\theta}=(\tilde{x},\tilde{z})$ to the system \eqref{eq:ms_inflated_avg_hybrid_sys} starting in ${K}_0+\delta\mathbb{B}$ there exists a solution $\bar{\theta}=(\bar{x},\bar{z})$ to system \eqref{eq:ms_avg_hybrid_sys} starting in $K_0$ such that the two solutions are $(T,\frac{\rho}{2})$-close. Such $\delta$ exists due to \cite[Proposition 6.34]{bookHDS}. Without loss of generality, we may take $\delta<1$ and $\rho<1$. 

Let $\mathcal{S}_{\mathcal{H}_2^{\text{ave}}}(K_0)$ be the set of maximal solutions to  (\ref{eq:ms_avg_hybrid_sys}) starting in $K_0$, and define the following sets:
\begin{align}
        \mathcal{R}_{T}(K_0)&:= \Big\{\xi=\bar{\theta}(t,j):\bar{\theta}\in\mathcal{S}_{\mathcal{H}_2^{\text{ave}}}({K}_0),\notag\\
        &~~~~~~~~~(t,j)\in\text{dom}(\bar{\theta}),~t+j\leq T\Big\}, \notag\\
        K_1 &:= \mathcal{R}_{T}(K_0) + \mathbb{B},~~~K:=K_1\cup G\left(K_1\cap D\right),\label{constructionsetK}
\end{align}
where ${K}$ is compact by \cite[Proposition 2]{TeelNesicAveraging}. Consider the functions $\psi_m,\psi_p:C_{\delta}\times\mathbb{R}^2_{\geq0}\to\mathbb{R}^{n_1}$ defined as follows:

\vspace{-0.2cm}
\begin{subequations}\label{auxpm}
\begin{align}
    \psi_m(\theta,\tau)&:= \frac{1}{2}\left( \partial_x \phi_1 \cdot u_1 - \partial_x u_1 \cdot \phi_1\right), \\
    \psi_p(\theta,\tau)&:=\frac{1}{2}\left( \partial_x \phi_1 \cdot u_1 + \partial_x u_1 \cdot \phi_1\right), 
\end{align}
\end{subequations}
where for simplicity we omitted the arguments in the right-hand side of \eqref{auxpm}. Using integration by parts and item (c) in Assumption \ref{asmp:A1}, we have  $\int_0^{T_2}\psi_p(\theta,\tau)\,d\tau_2= 0$, for all $(\theta,\tau_1)\in C_{\delta}\times\mathbb{R}_{\geq0}$. Let $\bar{h}:C_{\delta}\times\mathbb{R}_{\geq0}\to\mathbb{R}^{n_1}$ be defined as:
\begin{equation}\label{def:barh}
    \bar{{h}}(\theta,\tau_1):= \frac{1}{T_2}\int_0^{T_2}\Big(\phi_2(\theta,\tau)+\psi_m(\theta,\tau)\Big)\,d{\tau}_2,
\end{equation}
which satisfies $\bar{{f}}(\theta)= \frac{1}{T_1}\int_0^{T_1}\bar{{h}}(\theta,\tau_1)\,d\tau_1$. 

\vspace{0.1cm}
Below, Lemmas \ref{lemma0technical}-\ref{keytechnicallemma1} follow directly by the continuity and periodicity properties of the respective functions. For completness, the proofs can be found in the Appendix of the supplemental material, or in the Extended Manuscript \cite{abdelgalil2023lie}.

\begin{lemma}\label{lemma0technical}
Let $K\subset\mathbb{R}^n$ be a compact set. Then, there exist $L_{\bar{h}},L_{\bar{f}},M_{\bar{h}}>0$ such that $|\bar{h}(\theta,\tau_1)|\leq M_{\bar{h}}$ and $|\bar{h}(\theta,\tau_1)-\bar{h}(\theta',\tau_1)|\leq L_{\bar{h}}|\theta-\theta'|$, $|\bar{f}(\theta)-\bar{f}(\theta')|\leq L_{\bar{f}}|\theta-\theta'|$, for all $\theta,\theta'\in K\cap C_{\delta}$ and all $\tau_1\in \mathbb{R}_{\geq0}$. \QEDB 
\end{lemma}
%


\vspace{0.1cm}
\noindent Next, for each $(\theta,\tau)\in C_{\delta}\times\mathbb{R}^2_+$, we define
\begin{align*}
    {u}_{2}(\theta,\tau)&:= \int_0^{\tau_2} \Big(\phi_2(\theta,\tau_1,s_2)+\psi_m(\theta,\tau_1,s_2)\\
    &~~~~~~~~~~~~~+\psi_p(\theta,\tau_1,s_2) -\bar{{h}}(\theta,\tau_1,s_2)\Big)\,ds_2.
\end{align*}
as well as the auxiliary functions
\begin{subequations}\label{eq:defsigmas}
\begin{align}
    {\sigma}_{1} (\theta,\tau) &:= {u}_{1}(\theta,\tau),\label{constructionsigma1}\\
    {\sigma}_{2}(\theta,\tau)&:= {u}_{2}(\theta,\tau)-\int_{0}^{\tau_2}\partial_{\tau_1}{u}_1(\theta,\tau_1,s_2) \,ds_2\notag\\
    &~~~~~~~~~~~~~~- \partial_x{u}_1(\theta,\tau)\cdot {u}_{1}(\theta,\tau),\\
    \bar{{\sigma}}(\theta,\tau_1)&:= \int_0^{\tau_1}\Big(\bar{{h}}(\theta,s_1)-\bar{{f}}(\theta)\,\Big)\,ds_1.\label{eq:barsigma}
\end{align}
\end{subequations}
The functions $(\sigma_k,\bar{\sigma})$, $k\in\{1,2\}$, will play an important role in our analysis. We state some of their regularity properties.

\vspace{0.1cm}
\begin{lemma}\label{keylemma1periodic}
Let $K\subset\mathbb{R}^n$ be a compact set. Then, there exist $L_k,\bar{L}>0$ such that for $k\in\{1,2\}$:
\begin{subequations}\label{smoothnessbounds2}
\begin{align}
&|\sigma_k(\theta,\tau)|\leq L_k,~~~~|\bar{\sigma}(\theta,\tau_1)|\leq \bar{L},\label{sigmakprop1}\\
&|\sigma_k(\theta,\tau)-\sigma_k(\hat{\theta},\hat{\tau})|\leq L_k(|\theta-\hat{\theta}|+|\tau-\hat{\tau}|),\\
&|\bar{\sigma}(\theta,\tau_1)-\bar{\sigma}(\hat{\theta},\hat{\tau}_1)|\leq L(|\theta-\hat{\theta}|+|\tau_1-\hat{\tau}_1|),
\end{align}
\end{subequations}
$\forall~\theta,\hat{\theta}\in K\cap C_{\delta}$ $\tau,\hat{\tau}\in\mathbb{R}_{\geq0}^2$, and $\tau_1,\hat{\tau}_1\in\mathbb{R}_{\geq0}$.\QEDB 
\end{lemma}

%

\vspace{0.05cm}
\noindent The following technical Lemma will be key for our results. The proof of the first items leverages the Lipschitz extension Lemma of \cite[Lemma 2]{TeelNesicAveraging}. 

\vspace{0.05cm}
\begin{lemma}\label{keytechnicallemma1}
Let $K\subset\mathbb{R}^n$ be a compact set, and for $k\in\{1,2\}$ let $\sigma_k,\bar{\sigma}$ be given by \eqref{eq:defsigmas}, and $L_k,\bar{L}>0$ be given by Lemma \ref{keylemma1periodic}. Consider the closed sets
\begin{align}\label{eq:restrictedsets}
        C|_{K,\delta}=\left(K\cap C_{\delta}\right)\times\mathbb{R}^2_+,~~~~D|_{K}= \left(K\cap\bar{D}\right)\times\mathbb{R}^2_+.
\end{align}
Then, there exist functions $\tilde{\sigma}_k:\mathbb{R}^{n}\times\mathbb{R}^2_{+}\to\mathbb{R}^{n_1}$, $\tilde{\sigma}:\mathbb{R}^{n}\times\mathbb{R}_{\geq0}^2\to\mathbb{R}^{n_1}$ such that the following holds:
\begin{enumerate}[(a)]
\item For all $(\theta,\tau)\in {C}|_{K,\delta}$, we have $\widetilde{{\sigma}}_k(\theta,\tau)={\sigma}_k(\theta,\tau)$ and $\widetilde{{\sigma}}(\theta,\tau_1)=\bar{{\sigma}}(\theta,\tau_1)$.
\item  For all $(\theta,\tau)\in\mathbb{R}^{n+2}$, we have $\lvert \widetilde{{\sigma}}_k(\theta,\tau)\lvert\leq L_k$, and $\lvert \widetilde{{\sigma}}(\theta,\tau_1)\lvert\leq \bar{L}$. 
%
%
\item For all $(\theta,\tau)\in\mathbb{R}^{n+2}$ and $(\theta',\tau')\in\mathbb{R}^{n+2}$:
\begin{align*}
&\lvert \widetilde{{\sigma}}_k(\theta,\tau) - \widetilde{{\sigma}}_k(\theta',\tau')\lvert\leq L_k\,\big(\lvert \theta-\theta'\lvert+\lvert \tau-\tau'\lvert\big),\\
&\lvert \widetilde{{\sigma}}(\theta,\tau_1)- \widetilde{{\sigma}}(\theta',\tau'_1)\lvert \leq \sqrt{n}\bar{L}\,\big(\lvert \theta-\theta'\lvert+\lvert \tau_1-\tau'_1\lvert\big).
 \end{align*}
\item For all $(\theta,\tau)\in {C}|_{K,\delta}$:
\begin{align}
\alpha_1(\theta,\tau):&=\phi_1(\theta,\tau)-\partial_{\tau_2}\widetilde{{\sigma}}_1(\theta,\tau)=0,\label{importantPropLemma2}\\
\alpha_2(\theta,\tau):&=\phi_2(\theta,\tau)- \partial_{\tau_2}\widetilde{{\sigma}}_2(\theta,\tau)\notag\\
&~~~~~~~-\partial_{x}\widetilde{{\sigma}}_1(\theta,\tau)\,\phi_1(\theta,\tau)-\partial_{\tau_1}\widetilde{{\sigma}}_1(\theta,\tau)\notag\\
&=\bar{h}(\theta,\tau_1).\label{important2PropLemma2}
\end{align}
\item There exists $M_3>0$ such that for all $(\theta,\tau)\in {C}|_{K,\delta}$, $\alpha_3(\theta,\tau)\subset M_3\mathbb{B}$, where
\begin{align}
\nonumber \alpha_3(\theta,\tau)
:=\big\{&X\,\phi_1(\theta,\tau)+ Z \, \zeta+\partial_x\widetilde{{\sigma}}_1(\theta,\tau)\,\phi_2(\theta,\tau)\\\nonumber 
 &+\partial_{\tau_1}\widetilde{{\sigma}}_2(\theta,\tau):\,\zeta\in\Phi(z), Z\in \partial_z\widetilde{\sigma}_1(\theta,\tau),\\\label{eq:alpha_3_defn}& X\in\partial_x\widetilde{\sigma}_2(\theta,\tau)\big\},
\end{align}
\item There exists $M_4>0$ such that for all $(\theta,\tau)\in {C}|_{K,\delta}$, $\alpha_4(\theta,\tau)\subset M_4\mathbb{B}$, where
\begin{align}
\nonumber\alpha_4(\theta,\tau):=\big\{&X\,\phi_2(\theta,\tau)+Z\,\zeta: \,\zeta\in\Phi(z),\\ \label{eq:alpha_4_defn} &Z\in \partial_z\widetilde{\sigma}_1(\theta,\tau),X\in\partial_x\widetilde{\sigma}_2(\theta,\tau)\big\}.
\end{align}
\end{enumerate}
\end{lemma}


\vspace{0.2cm}
\noindent Continuing with the proof of Theorem \ref{thm:closeness_of_trajs}, let $\varepsilon_1^*:=\delta(\sqrt{n}\bar{L}M_s+M_3+ M_4+0.5\bar{L}_{\bar{h}})^{-1}$, $\varepsilon_2^*:=\delta(2\bar{L}(0.5\sqrt{n}+1))^{-1}$, $\varepsilon_3^*:=\delta(
2(L_1+L_2)+1)^{-1}$, $\varepsilon_4^*:=\rho
(2(L_1+L_2+(0.5\sqrt{n}+1)\bar{L}))^{-1}$, and $\varepsilon^*:=\min\left\{\varepsilon_0,\varepsilon_1,\varepsilon_2,\varepsilon_3,\varepsilon_4\right\}$, which satisfies $\varepsilon^*\in(0,1)$ due to the definition of $\varepsilon_3$. Let $\varepsilon\in(0,\varepsilon^*)$ and consider the restricted HDS:

\vspace{-0.2cm}
\begin{subequations}\label{eq:ms_orig_hybrid_sys_restricted}
\begin{align}
(\theta,\tau)&\in C|_{K},~~~\left\{\begin{array}{l}
~\dot{\theta}~\in~ F_{\varepsilon}(\theta,\tau_1,\tau_2),\\
\dot{\tau}_1~=~\varepsilon^{-1},\\
\dot{\tau}_2~=~\varepsilon^{-2},
\end{array}\right.\label{flowmaprestricted}\\
(\theta,\tau)&\in D|_{K},~~~\left\{\begin{array}{l}
\theta^+~\in~ G(\theta)\cap K,\\
\tau^+_1~=~\tau_1,\\
\tau^+_2~=~\tau_2,
\end{array}\right.\label{jumpmaprestricted}
\end{align}
\end{subequations}
with $\theta=(x,z)$, $\tau=(\tau_1,\tau_2)$, and flow set $C|_{K}$ and jump set $D|_{K}$ defined in \eqref{eq:restrictedsets}, with $K$ given by \eqref{constructionsetK}. We divide the rest of proof into four main steps:

\vspace{0.1cm}
\noindent 
\textsl{Step 1: Construction and Properties of Auxiliary Functions:} Using the compact set $\tilde{K}:=K+\frac{\delta}{2}\mathbb{B}$ and the functions $\sigma_k,\bar{\sigma}$ given by \eqref{eq:defsigmas}, let Lemmas \ref{lemma0technical}-\ref{keylemma1periodic} generate   $L_{\bar{h}},L_{\bar{f}},M_{\bar{h}},L_k,\bar{L}>0$. Using these constants, the set $\tilde{K}$, and the functions ($\sigma_k,\bar{\sigma}$), let Lemma \ref{keytechnicallemma1} generate the functions $(\tilde{\sigma}_k,\tilde{\sigma})$. Then, for each $(\theta,\tau)\in\mathbb{R}^n\times\mathbb{R}^2_{+}$ and each $\varepsilon\in(0,\varepsilon^*)$, we define
\begin{equation}\label{constructiontheta}
\hat{\theta}:=\theta-\gamma_{\varepsilon}(\theta,\tau),~~~
\tilde{\theta}:=\hat{\theta}-\lambda_{\varepsilon}(\hat{\theta},\tau_1),
\end{equation}
where
\begin{align}
\gamma_{\varepsilon}(\theta,\tau):&=\left(\!\!\begin{array}{c}
\sum_{k=1}^2\varepsilon^k\tilde{\sigma}_k(\theta,\tau)\\
0
\end{array}\!\!\!\right),~\lambda_{\varepsilon}(\hat{\theta},\tau_1):=\left(\!\!\!\begin{array}{c}
\varepsilon\tilde{\sigma}(\hat{\theta},\tau_1)\\
0
\end{array}\!\!\!\right).\label{deflambdae}
\end{align}
Using the inequalities in item (b) of Lemma \ref{keytechnicallemma1}, we have 
\begin{align}\label{auxdefgammaepsilon}
\left|\gamma_{\varepsilon}(\theta,\tau)\right|&=\left|\sum_{k=1}^2\varepsilon^k \widetilde{{\sigma}}_k(\theta,\tau)\right|\leq \varepsilon(L_1+\varepsilon L_2),
\end{align}
for all $(\theta,\tau)\in \mathbb{R}^n\times\mathbb{R}^2_+ $. Similarly, due to items (b)-(c) of Lemma \ref{keytechnicallemma1}, the definition of $\hat{\theta}$ in \eqref{constructiontheta}, and inequality \eqref{auxdefgammaepsilon},
\begin{align}\label{auxdeflambdaepsilon}
\left|\lambda_{\varepsilon}(\hat{\theta},\tau_1)\right|&\leq \varepsilon\left|\tilde{\sigma}(\hat{\theta},\tau_1)-\tilde{\sigma}\left(\theta,\tau_1\right)\right|+\varepsilon\Big|\tilde{\sigma}\left(\theta,\tau_1\right)\Big|\notag\\
&\leq \varepsilon\sqrt{n}\bar{L}\big|\hat{\theta}-\theta\big|+\varepsilon|\tilde{\sigma}(\theta,\tau_1)|\notag\\
&\leq \varepsilon\sqrt{n}\bar{L}\big|\gamma_{\varepsilon}(\theta,\tau)\big|+\varepsilon|\tilde{\sigma}(\theta,\tau_1)|,
\end{align}
for all $(\theta,\tau)\in\mathbb{R}^n\times\mathbb{R}_{\geq0}^2$. Also, note that due to item (a) of Lemma \ref{keytechnicallemma1}, the functions $\gamma_{\varepsilon}$ and  $\lambda_{\varepsilon}$ satisfy 
$\gamma_{\varepsilon}(\theta,\tau)=\left(\sum_{k=1}^2\varepsilon^k\sigma_k(\theta,\tau),0\right)$ and $\lambda_{\varepsilon}(\hat{\theta},\tau_1)=\left(\varepsilon\bar{\sigma}(\hat{\theta},\tau_1),0\right)$ for all $(\theta,\tau),(\hat{\theta},\tau)\in \left(\tilde{K}\cap C_{\delta}\right)\times\mathbb{R}^2_+$, with $(\sigma_k,\bar{\sigma})$ given by \eqref{smoothnessbounds2}.

\vspace{0.2cm}\noindent 
\textsl{Step 2: Construction of First Auxiliary Solution:} Let $(\theta,\tau)$ be a solution to \eqref{eq:ms_orig_hybrid_sys_restricted} with $\theta(0,0)\in K_0$. By construction, we have
\begin{equation}\label{inclusionthetak}
\theta(t,j)\in K\subset\tilde{K},~~\forall~(t,j)\in\text{dom}(\theta,\tau).
\end{equation}
Using \eqref{sigmakprop1}, item (b) of Lemma \ref{keytechnicallemma1}, \eqref{auxdefgammaepsilon}, and the choice of $\varepsilon^*$
\begin{equation}\label{deltabound001}
|\gamma_{\varepsilon}(\theta(t,j),\tau(t,j))|\leq \frac{\delta}{2},~~~\forall~(t,j)\in\text{dom}(\theta,\tau).
\end{equation}
For each $(t,j)\in\text{dom}(\theta,\tau)$, let $\hat{\theta}$ be defined via \eqref{constructiontheta}. It follows that $\hat{\theta}$ is a hybrid arc, and due to \eqref{inclusionthetak} and \eqref{deltabound001} it satisfies
\begin{equation}\label{boundthetahatcompact}
\hat{\theta}(t,j)\in K+\frac{\delta}{2}\mathbb{B}= \tilde{K},~~~~\forall~(t,j)\in\text{dom}(\hat{\theta},\tau).
\end{equation}
Thus, using \eqref{sigmakprop1}, item (b) of Lemma \ref{keytechnicallemma1}, and \eqref{auxdeflambdaepsilon}, we get
\begin{equation}\label{boundlambdalongtrajectories}
|\lambda_{\varepsilon}(\hat{\theta}(t,j),\tau_1(t,j))|\leq \varepsilon\bar{L} \left({\frac{1}{2}}\sqrt{n}+1\right)\leq{\frac{\delta}{2}},
\end{equation}
for all $(t,j)\in\text{dom}(\hat{\theta},\tau)$. 

\vspace{0.1cm}
Next, for each $(t,j)\in\text{dom}(\hat{\theta},\tau)$ such that $(t,j+1)\in\text{dom}(\hat{\theta},\tau)$, we have that: 
\begin{align*}
\theta(t,j)=\hat{\theta}(t,j)+\gamma_{\varepsilon}(\theta(t,j),\tau(t,j))\in D\cap K\subset D\cap \tilde{K}. 
\end{align*}
Using \eqref{deltabound001} and the construction in \eqref{inflatedjumpset}, we conclude that $\hat{\theta}(t,j)\in D_{\delta/2}$. Thus, using  \eqref{jumpmaprestricted} and  \eqref{deltabound001}, we obtain:
\begin{align}
\hat{\theta}(t,j+1)&=\theta(t,j+1)-\gamma_{\varepsilon}(\theta(t,j+1),\tau(t,j+1))\notag\\
&\in G\left(\theta(t,j)\cap D\right)\cap K+{\frac{\delta}{2}}\mathbb{B} \notag\\
    &\subset G\left(\theta(t,j)\cap D\right)+{\frac{\delta}{2}}\mathbb{B}\notag\\
    &= G\left(\left(\hat{\theta}(t,j)+\gamma_{\varepsilon}(\theta(t,j),\tau(t,j))\right)\cap D\right) + {\frac{\delta}{2}}\mathbb{B}\notag\\
     &\subset G\left(\left(\hat{\theta}(t,j)+{\frac{\delta}{2}}\mathbb{B}\right)\cap D\right) + {\frac{\delta}{2}}\mathbb{B}.
    \label{inflationjumps1}
\end{align}
Similarly, from \eqref{constructiontheta}, for all $j\in\mathbb{Z}_{\geq0}$ such that $I_j=\{t:\,(t,j)\in\text{dom}\,(\hat{\theta})\}$ has a nonempty interior, we have:
\begin{equation}\label{flowsetbelongs1}
\theta(t,j)=\hat{\theta}(t,j)+\gamma_{\varepsilon}(\theta(t,j),\tau(t,j))\in C\cap K\subset C_{\delta}\cap \tilde{K},
\end{equation}
for all $t\in I_j$. Therefore, using the construction of the inflated flow set \eqref{inflatedflowset}, and the bound \eqref{deltabound001}, we have that $\hat{\theta}(t,j)\in C_{{\delta/2}}$. Since $\gamma_{\varepsilon}$ is Lipschitz continuous due to Lemma \ref{keytechnicallemma1},  $\hat{\theta}(\cdot,j)$ is locally absolutely continuous, and it satisfies 
\begin{align}
\dot{\hat{\theta}}(t,j)&=\dot{\theta}(t,j)-\dot{\overbrace{\gamma_{\varepsilon}(\theta(t,j),\tau(t,j))}}\notag\\
&= F_{\varepsilon}(\theta,\tau)-\dot{\overbrace{\gamma_{\varepsilon}(\theta(t,j),\tau(t,j))}},\label{expandingflows}
\end{align}
for almost all $t\in I_j$. To compute \eqref{expandingflows}, note that
\begin{align*}
\varepsilon\dot{\tilde{\sigma}}_1&\in\varepsilon\partial_{\theta} \tilde{\sigma}_1\dot{\theta}+\varepsilon\partial_{\tau} \tilde{\sigma}_1\dot{\tau}\\
&=\varepsilon\partial_{x} \tilde{\sigma}_1\dot{x}+\varepsilon\partial_{z} \tilde{\sigma}_1\dot{z}+\varepsilon\partial_{\tau_1} \tilde{\sigma}_1\dot{\tau}_1+\varepsilon\partial_{\tau_2} \tilde{\sigma}_1\dot{\tau}_2\\
&\subset \partial_{x} \tilde{\sigma}_1\sum_{k=1}^2\varepsilon^{k-1}{\phi}_{k}+\partial_{\tau_1} \tilde{\sigma}_1+\frac{1}{\varepsilon}\partial_{\tau_2} \tilde{\sigma}_1+\varepsilon\partial_{z} \tilde{\sigma}_1\Phi,\\
\varepsilon^2\dot{\tilde{\sigma}}_2&\in\varepsilon^2\partial_{\theta} \tilde{\sigma}_2\dot{\theta}+\varepsilon^2\partial_{\tau} \tilde{\sigma}_2\dot{\tau}\\
&=\varepsilon^2\partial_{x} \tilde{\sigma}_2\dot{x}+\varepsilon^2\partial_{z} \tilde{\sigma}_2\dot{z}+\varepsilon^2\partial_{\tau_1} \tilde{\sigma}_2\dot{\tau}_1+\varepsilon^2\partial_{\tau_2} \tilde{\sigma}_2\dot{\tau}_2\\
&\subset \partial_{x} \tilde{\sigma}_2\sum_{k=1}^2\varepsilon^{k}{\phi}_{k}+\varepsilon\partial_{\tau_1} \tilde{\sigma}_2+\partial_{\tau_2} \tilde{\sigma}_2+\varepsilon^2\partial_{z} \tilde{\sigma}_2\Phi,
\end{align*}
where we used \eqref{flowmapstructure}, the structure of \eqref{vectofield0}, and \eqref{eq:ms_orig_hybrid_sys_restricted}. Therefore, using the definition of $\gamma_{\varepsilon}$, we have that \eqref{expandingflows} can be written as 
\begin{gather*}
    \dot{\hat{\theta}}(t,j)\in \left(\begin{array}{c}\sum_{k=1}^4\varepsilon^{k-2}{\alpha}_k(\theta(t,j),\tau(t,j))\\
    \Phi(z(t,j))
    \end{array}\right),
\end{gather*}
for almost all $t\in I_j$. Using the last containment in \eqref{flowsetbelongs1}, and items (a), (d), (e), and (f) of Lemma \ref{keytechnicallemma1}, we obtain 
\begin{align*}
&\alpha_1(\theta(t,j),\tau(t,j))=0,\\
&\alpha_2(\theta(t,j),\tau(t,j))=\bar{{h}}(\theta(t,j),\tau_1(t,j)),\\
&\alpha_k(\theta(t,j),\tau(t,j))\subset M_k\mathbb{B},~~~k\in\{3,4\}.
\end{align*}
for all $t\in I_j$. Thus, $\forall~\varepsilon\in(0,\varepsilon^*)$ and
\begin{align}\label{dynamicshatxie}
    \dot{\hat{\theta}}(t,j)\in \left(\begin{array}{c} 
    \bar{{h}}(\theta(t,j),\tau_1(t,j))\\
    \Phi(z(t,j))
    \end{array}\right)+
    \chi_{\varepsilon}(\theta(t,j),\tau(t,j)),
\end{align}
for almost all $t\in I_j$, where
\begin{equation*}
\chi_{\varepsilon}(\theta,\tau):=\left(\begin{array}{c}
\sum_{k=1}^2\varepsilon^k\alpha_{k+2}(\theta(t,j),\tau(t,j))\\
0
\end{array}\right).
\end{equation*}
Finally, using the equality in \eqref{flowsetbelongs1}, we can write \eqref{dynamicshatxie} as: 
\begin{align}
\dot{\hat{\theta}}(t,j)&\in\left(\begin{array}{c} 
    \bar{{h}}\left(\hat{\theta}(t,j)+\gamma_{\varepsilon}(\theta(t,j),\tau(t,j)),\tau_1(t,j)\right)\\
    \Phi(z(t,j))
    \end{array}\right)\notag\\
    & ~~~~~~~~~~~~~~~~~~~~~~~~~~~~~+\chi_{\varepsilon}(\theta(t,j),\tau(t,j)),\label{dothatfirst}
\end{align}
which holds for almost all $t\in I_j$.

\vspace{0.1cm}
\noindent 
\textsl{Step 3: Construction of Second Auxiliary Solution:} Using the solution $(\theta,\tau)$ to the restricted HDS \eqref{eq:ms_orig_hybrid_sys_restricted} considered in Step 2, as well as the hybrid arc $\hat{\theta}$, we now define a new auxiliary solution. In particular, for each $(t,j)\in\text{dom}(\theta,\tau)$, using \eqref{constructiontheta} we define:
\begin{equation}\label{constructiontildetheta}
\tilde{\theta}(t,j):=\hat{\theta}(t,j)-\lambda_{\varepsilon}\left(\hat{\theta}(t,j),\tau_1(t,j)\right).
\end{equation}
By construction, $\tilde{\theta}$ is a hybrid arc. Since, by  combining equation \eqref{constructiontheta} and inequalities  \eqref{deltabound001} and  \eqref{boundlambdalongtrajectories} we have that
\begin{equation}\label{delta2bound}
|\gamma_{\varepsilon}(\theta(t,j),\tau(t,j))+\lambda_{\varepsilon}(\hat{\theta}(t,j),\tau_1(t,j))|\leq {\delta},
\end{equation}
for all $(t,j)\in\text{dom}(\theta,\tau)$, and for all $(t,j)\in\text{dom}(\hat{\theta},\tau)$, it follows that the hybrid arc \eqref{constructiontildetheta} satisfies $\tilde{\theta}(0,0)\in K_0+\delta\mathbb{B}$. 

\vspace{0.1cm}
Now, for each $(t,j)\in\text{dom}(\tilde{\theta})$ such that $(t,j+1)\in\text{dom}(\tilde{\theta})$ it satisfies $\theta(t,j)=\tilde{\theta}(t,j)+\gamma_{\varepsilon}(\theta(t,j),\tau(t,j))+\lambda_{\varepsilon}(\hat{\theta}(t,j),\tau_1(t,j))\in D\cap K\subset D\cap\tilde{K}$, where we used \eqref{constructiontheta} and \eqref{constructiontildetheta}. Therefore, using the construction of the inflated jump set \eqref{inflatedjumpset}, we conclude that $\tilde{\theta}(t,j)\in D_{{\delta}}$, and using \eqref{boundlambdalongtrajectories}-\eqref{inflationjumps1} we have:
\begin{align*}
\tilde{\theta}(t,j+1)&=\hat{\theta}(t,j+1)-\lambda_{\varepsilon}(\hat{\theta}(t,j+1),\tau_1(t,j+1))\\
&\in G\left(\left(\hat{\theta}(t,j)+{\frac{\delta}{2}}\mathbb{B}\right)\cap D\right)+{\delta}\mathbb{B}\\
&=G\Big(\Big(\tilde{\theta}(t,j)+\lambda_{\varepsilon}(\hat{\theta}(t,j),\tau_1(t,j))+{\frac{\delta}{2}}\mathbb{B}\Big)\cap D\Big)\\
&~~~+{\delta}\mathbb{B} \subset G_{{\delta}}(\tilde{\theta}(t,j)),
\end{align*}
where we also used \eqref{constructiontildetheta} and the definition of the inflated jump map $G_{{\delta}}$ in \eqref{inflatiojumpmap}. On the other hand, for each $j\in\mathbb{Z}_{\geq0}$ such that $I_j=\{t:\,(t,j)\in\text{dom}\,(\tilde{\theta})\}$ has a nonempty interior, we must have:
\begin{align}
\theta(t,j)&=\tilde{\theta}(t,j)+\gamma_{\varepsilon}(\theta(t,j),\tau(t,j))+\lambda_{\varepsilon}(\hat{\theta}(t,j),\tau_1(t,j))\notag\\
&=\hat{\theta}(t,j)+\lambda_{\varepsilon}(\hat{\theta}(t,j),\tau_1(t,j))\in C\cap K,\label{thetainclusion0k}
\end{align}
and therefore, due to \eqref{boundlambdalongtrajectories}, the inclusion in \eqref{thetainclusion0k}, and the definition of inflated flow set \eqref{inflatedflowset}: 
\begin{equation}\label{thetahatinclusion}
\hat{\theta}(t,j)\in C_{{\delta/2}}\cap K\subset C_{{\delta/2}}\cap \tilde{K}\subset C_{\delta}\cap\tilde{K},~~~\forall~t\in I_j.
\end{equation}
Since $\tilde{\sigma}$ is Lipschitz continuous in all arguments due to Lemma \ref{keytechnicallemma1}, it follows that $\tilde{\theta}(\cdot,j)$ is locally absolutely continuous and satisfies:
\begin{equation*}
\dot{\tilde{\theta}}(t,j)=\dot{\hat{\theta}}(t,j)-\frac{d\lambda_{\varepsilon}(\hat{\theta}(t,j),\tau_1(t,j))}{dt},~~\text{for~a.a.}~~t\in I_j.
\end{equation*}
Using the definition of $\lambda_{\varepsilon}$ in \eqref{deflambdae}, as well as \eqref{thetahatinclusion}, item (a) of Lemma \ref{keytechnicallemma1}, and \eqref{eq:barsigma}, we obtain:
\begin{align*}
\dot{\lambda}_{\varepsilon}&\in\left(\varepsilon \partial_{\hat{\theta}}\tilde{\sigma}(\hat{\theta},\tau_1)\dot{\hat{\theta}}+\varepsilon \partial_{\tau_1}\tilde{\sigma}(\hat{\theta},\tau_1)\dot{\tau}_1,0\right)\\
&=\left(\varepsilon \partial_{\hat{\theta}}\tilde{\sigma}(\hat{\theta},\tau_1)\dot{\hat{\theta}}+ \partial_{\tau_1}\tilde{\sigma}(\hat{\theta},\tau_1),0\right)\\
&=\left(\varepsilon \partial_{\hat{\theta}}\tilde{\sigma}(\hat{\theta},\tau_1)\dot{\hat{\theta}}+ \bar{h}(\hat{\theta},\tau_1)-\bar{f}(\hat{\theta}),0\right).
\end{align*}
Therefore, using \eqref{dothatfirst} we have that for almost all $t\in I_j$:
\begin{align}\label{inflationtildetheta}
\dot{\tilde{\theta}}(t,j)&\in  \left(\!\!\begin{array}{c} 
    \bar{f}\big(\hat{\theta}(t,j)\big)\\
    \Phi\left(z(t,j)\right)
    \end{array}\!\!\right)+\left(\!\!\!\begin{array}{c}
    p\left(\theta(t,j),\hat{\theta}(t,j),\tau(t,j)\right)\\
    0
    \end{array}\!\!\!\right),
\end{align}
where the last term can be written in compact form as
\begin{align*}
p(\theta,\hat{\theta},\tau)&=\varepsilon \partial_{\hat{\theta}}\tilde{\sigma}(\hat{\theta},\tau_1)\dot{\hat{\theta}}+
\sum_{k=1}^2\varepsilon^k\alpha_{k+2}(\theta,\tau)\\
&~~~~~~~+\bar{{h}}\left(\hat{\theta}+\gamma_{\varepsilon}(\theta,\tau),\tau_1\right)-\bar{h}(\hat{\theta},\tau_1).
\end{align*}
Using \eqref{constructiontildetheta}, \eqref{thetainclusion0k}, \eqref{thetahatinclusion} and Lemma \ref{lemma0technical}, we have:
\begin{subequations}\label{inequalitiesusefulfinal}
\begin{align}
    \left\lvert \bar{h}\left(\hat{\theta}+\gamma_{\varepsilon}(\theta,\tau),\tau_1\right)-\bar{h}\left(\hat{\theta},\tau_1\right)\right\lvert &\leq L_{\bar{h}}\lvert  \gamma_{\varepsilon}(\theta,\tau)\lvert,\label{inequalitiesusefulfinal1}\\
    \left\lvert \bar{{h}}(\theta,\tau_1)\right\lvert &\leq M_{\bar{h}}.\label{usefulboundhbar}
\end{align}
\end{subequations}
for all $t\in I_j$ and all $\varepsilon\in(0,\varepsilon^*)$, where we omitted the time arguments to simplify notation. Using $\hat{\theta}=(\hat{x},z)$, we have that 
\begin{align*}
\partial_{\hat{\theta}}\tilde{\sigma}(\hat{\theta},\tau_1)\dot{\hat{\theta}}&=\partial_{\hat{x}}\tilde{\sigma}(\hat{\theta},\tau_1)\dot{\hat{x}}+\partial_{z}\tilde{\sigma}(\hat{\theta},\tau_1)\dot{z}\\
&\subset \partial_{\hat{x}}\tilde{\sigma}(\hat{\theta},\tau_1)\left(\bar{h}(\theta,\tau_1)+\sum_{k=1}^2\varepsilon^k\alpha_{k+2}(\theta,\tau)\right)\\
&~~~~~~+\partial_{z}\tilde{\sigma}(\hat{\theta},\tau_1)\Phi(z),
\end{align*}
where the inclusion follows from \eqref{dothatfirst}, and $\alpha_i$ are defined in \eqref{importantPropLemma2}, \eqref{important2PropLemma2}, \eqref{eq:alpha_3_defn}, and \eqref{eq:alpha_4_defn}. Since $\tilde{\sigma}$ is globally Lipschitz continuous by Lemma \ref{keytechnicallemma1}, it follows that $ \partial_{\hat{x}}\tilde{\sigma}\subset \bar{L}\mathbb{B}, \,\partial_{\hat{z}}\tilde{\sigma}\subset \bar{L}\mathbb{B}$. Using  \eqref{thetainclusion0k}, \eqref{thetahatinclusion}, items (c), (e) and (f) of Lemma \ref{keytechnicallemma1}, as well as \eqref{usefulboundhbar}, we obtain for all $t\in I_j$ and all $\varepsilon\in(0,\varepsilon^*)$:
\begin{align}\label{boundperturbationfinal}
\partial_{\hat{\theta}}\tilde{\sigma}(\hat{\theta}(t,j),\tau_1(t,j))\dot{\hat{\theta}}(t,j)
&\subset\sqrt{n}\bar{L}M_s \mathbb{B},
\end{align}
where $M_s:=M_{\bar{h}}+M_3+M_4+M_{z}$, the constant $M_{z}$ comes from the proof of items (e)-(f) in Lemma \ref{keytechnicallemma1}, and where we used the fact that $\varepsilon<1$. Therefore, using the above bounds, as well as \eqref{deltabound001} to bound \eqref{inequalitiesusefulfinal1}, we conclude that for all $t\in I_j$ and all $\varepsilon\in(0,\varepsilon^*)$: $p(\theta(t,j),\hat{\theta}(t,j),\tau(t,j))\subset\varepsilon\widetilde{M}\mathbb{B}\subset {\delta} \mathbb{B}$, where $\widetilde{M} = \sqrt{n}\bar{L}M_s+M_3+ M_4+{\frac{\delta\bar{L}_{\bar{h}}}{2}}$, and we used the fact that $\delta<1$ and the choice of $\varepsilon$. Using \eqref{boundlambdalongtrajectories}, \eqref{constructiontildetheta}, and \eqref{inflationtildetheta}, we conclude that for almost all $t\in I_j$:
\begin{align*}
\dot{\tilde{\theta}}(t,j)&\in \left(\begin{array}{c}
\bar{f}\left(\tilde{\theta}(t,j)+{\frac{\delta}{2}}\mathbb{B}\right)\\
\Phi(z(t,j))
\end{array}\right)+{\delta}\mathbb{B}\subset F_{\delta}(\tilde{\theta}(t,j)),
\end{align*}
where $F_{\delta}$ is the inflated average flow map in \eqref{eq:ms_inflated_avg_hybrid_sys}. It follows that $\tilde{\theta}$ is $(T,\rho/2)$-close to some solution $\bar{\theta}$ of the average system \eqref{eq:ms_avg_hybrid_sys}, with $\bar{\theta}(0,0)\in K_0$. Using \eqref{constructiontheta}, \eqref{constructiontildetheta}, and the fact that \eqref{delta2bound} also holds with $\delta=\frac{\rho}{2}$ due to the fact that $\varepsilon<\varepsilon_4$, we conclude that $\theta$ is $(T,{\rho})$-close to $\bar{\theta}$.

\vspace{0.15cm}
\textsl{Step 4: Removing $K$ from the HDS \eqref{eq:ms_orig_hybrid_sys_restricted}.} We now study the properties of the solutions to the unrestricted system \eqref{eq:ms_orig_hybrid_sys} from $K_0$, based on the properties of the solutions of the restricted system  \eqref{eq:ms_orig_hybrid_sys_restricted} initialized also in $K_0$. Let $(\theta,\tau)$ be a solution of \eqref{eq:ms_orig_hybrid_sys} starting in $K_0$. We consider two scenarios:
\begin{enumerate}[(a)]
\item For all $(t,j)\in\text{dom}(\theta,\tau)$ such that $t+j\leq T$ we have $\theta(t,j)\in K$. Then, it follows that $\theta$ is $(T,\rho)$-close to $\bar{\theta}$.
\item There exists $(t,j)\in\text{dom}(\theta,\tau)$ such that $\theta(s,k)\in K$ for all $(s,k)\in\text{dom}(\theta,\tau)$ such that $s+k\leq t+j$ and either: 
\begin{enumerate}[(1)]
\item there exists sequence $\{\ell_i\}_{i\in\mathbb{Z}}$ satisfying $\lim_{i\to\infty}\ell_i=t$ such that $(\ell_i,j)\in\text{dom}(\theta,\tau)$ and $\theta(\ell_i,j)\notin K$ for each $i$, or else
\item $(t,j+1)\in\text{dom}(\theta,\tau)$ and $\theta(t,j+1)\notin K$.
\end{enumerate}
Then, we must have that the solution $(\theta,\tau)$ agrees with a solution to \eqref{eq:ms_orig_hybrid_sys_restricted} up to time $(t,j)$, which implies that $\theta(t,j)\in \mathcal{R}_T(K_0)+\rho\mathbb{B}$. But, since $\rho\in(0,1)$, and due to the construction of $K$ in \eqref{constructionsetK}, we have that $\mathcal{R}_T(K_0)+\rho\mathbb{B}$ is contained in the interior of $K$, so neither of the above two cases can occur. 
\end{enumerate}
Since item (b) cannot occur, we obtain the desired result. \hfill $\blacksquare$

\vspace{-0.2cm}
\subsection{Proof of Theorem \ref{thm:avg_UAS_implies_org_PUAS}} The proof proceeds in a similar manner to the proof of \thref{thm:closeness_of_trajs} and the proof of \cite[Theorem 2]{TeelNesicAveraging}. 

Let $\omega:\mathcal{B}_{\mathcal{A}}\rightarrow\mathbb{R}_{\geq 0}$ be a proper indicator for $\mathcal{A}$ on $\mathcal{B}_{\mathcal{A}}$, and let $\beta\in\mathcal{KL}$ such that each solution $\bar{\theta}$ of the averaged HDS \eqref{eq:ms_avg_hybrid_sys} starting in $\mathcal{B}_{\mathcal{A}}$ satisfies 
\begin{align*}
    \omega\left(\bar{\theta}(t,j)\right)&\leq \beta\left(\omega(\bar{\theta}(0,0)),t+j\right), & \forall&(t,j)\in\text{dom}\,\bar{\theta}.
\end{align*}
Let $K_0\subset\mathcal{B}_{\mathcal{A}}$ be compact, and let
\begin{align}
    K_1&= \left\{\theta\in\mathcal{B}_{\mathcal{A}}:\, \omega(\theta)\leq \beta\left(\max_{\tilde{\theta}\in K_0}\omega(\tilde{\theta}),0\right) + 1\right\},\\
    K&= K_1\cup G(K_1\cap D).
\end{align}
By construction, the continuity of $\omega,\beta$, and the OSC property of $G$, the set $K$ is compact and satisfies $K\subset\mathcal{B}_{\mathcal{A}}$. Let $\nu\in(0,1)$ and observe that, due to the robustness properties of well-posed HDS \cite[Lemma 7.20]{bookHDS}, there exists $\delta\in(0,\delta^*)$ such that all solutions $\tilde{\theta}$ to the inflated averaged HDS \eqref{eq:ms_inflated_avg_hybrid_sys} that start in $K_0+\delta\mathbb{B}$ satisfy for all $(t,j)\in\text{dom}\,\tilde{\theta}$:
\begin{align}\label{inflatedover3}
    \omega(\tilde{\theta}(t,j))&\leq \beta\left(\omega(\tilde{\theta}(0,0)),t+j\right) + \frac{\nu}{3}.
\end{align}
Without loss of generality we may assume that $\delta<1$, and we define $\tilde{K}=K+\frac{\delta}{2}\mathbb{B}$. Using $\tilde{K}$ we let Lemmas \ref{lemma0technical} and \ref{keylemma1periodic} generate the constants $L_{\bar{h}},L_{\bar{f}},M_{\bar{h}},\bar{L},L_k,L>0$ so that the bounds Lemma \ref{lemma0technical} and \eqref{smoothnessbounds2} hold for all $\theta,\theta'\in \tilde{K}\cap C_{\delta}$ and all $\tau_1,\hat{\tau}_1\in\mathbb{R}_{\geq0}$. Using these constants, we define $\varepsilon_i$ as in the proof of Theorem 1, for all $i\in\{1,2,3,4\}$. Since $\omega$ and $\beta$ are continuous, and $\beta(r,\cdot)$ converges to zero as the argument grows unbounded, there exists $\varepsilon_5^*\in\left(0,\frac{\delta}{2(L_1+L_2+\bar{L})}\right)$ such that for all $~\theta\in \tilde{K}$ and $\tilde{\theta}\in \tilde{K}+\varepsilon_5^*(L_1+L_2+\bar{L})\mathbb{B}$ satisfying $|\theta-\tilde{\theta}|\leq \varepsilon_5^*(L_1+L_2+\bar{L})$ and all $s\geq 0$:
\begin{align}\label{inequalitiesindicator}
    \omega(\theta)&\leq \omega(\tilde{\theta}) + \frac{\nu}{3},  &
    \beta(\omega(\tilde{\theta}),s)&\leq \beta(\omega({\theta}),s) + \frac{\nu}{3}.
\end{align}
Then, we let $\varepsilon^* = \min_{k=1,\ldots,5}\varepsilon_k$ and fix $\varepsilon\in(0,\varepsilon^*)$. As in the proof of Theorem \ref{thm:closeness_of_trajs}, we define the restricted HDS:
\begin{subequations}\label{eq:ms_avg_hds_restricted_stab1}
\begin{align}
(\theta,\tau)&\in C|_{K}  & &\begin{cases}    
~~~~\dot{\theta}&\in{F}_{\varepsilon}(\theta,\tau_1,\tau_2),\\
    \hphantom{\Big[}\begin{array}{c}\dot{\tau}_1\end{array}\hphantom{\Big]}&= \varepsilon^{-1},\\
    \hphantom{\Big[}\begin{array}{c}\dot{\tau}_2\end{array}\hphantom{\Big]}&= \varepsilon^{-2},
\end{cases}\\
(\theta,\tau)&\in D|_{K} & &\begin{cases} 
~~~~\theta^+&\in G(\theta)\cap K,\\
\hphantom{\Big[}\begin{array}{c}\tau_1^+\end{array}\hphantom{\Big]} &= \tau_1,\\
\hphantom{\Big[}\begin{array}{c}\tau_2^+\end{array}\hphantom{\Big]} &= \tau_2,
\end{cases}
\end{align}
\end{subequations}
and we let $(\theta,\tau)$ denote a solution to \eqref{eq:ms_avg_hds_restricted_stab1} starting in $K_0$. As in the proof of Theorem \ref{thm:closeness_of_trajs}, using \eqref{constructiontheta}, we define, for each $(t,j)\in\text{dom}(\theta,\tau)$, the arc $
\tilde{\theta}:=\theta-\gamma_{\varepsilon}(\theta,\tau)-\lambda_{\varepsilon}(\theta-\gamma_{\varepsilon}(\theta,\tau),\tau_1)$, which, as shown in the proof of Theorem \ref{thm:closeness_of_trajs}, is a hybrid arc that is a solution to the inflated HDS \eqref{eq:ms_inflated_avg_hybrid_sys}, and therefore satisfies \eqref{inflatedover3} for all $(t,j)\in\text{dom}\,\tilde{\theta}$. We can now use \eqref{inequalitiesindicator} to obtain:
\begin{align}
    \omega(\theta(t,j))&\leq \omega(\tilde{\theta}(t,j))+\frac{\nu}{3}\notag\leq \beta(\omega(\tilde{\theta}(0,0)),t+j) + \frac{2\nu}{3} \notag\\
    &\leq \beta(\omega({\theta}(0,0)),t+j) + \nu,\label{thetaboundkl}
\end{align}
for all $(t,j)\in\text{dom}(\theta)$. Since  $\beta\in\mathcal{K}\mathcal{L}$, $\theta$ remains in the set  $K_\nu= \big\{\theta\in\mathcal{B}_{\mathcal{A}}:\, \omega(\theta)\leq \beta\big(\max_{\tilde{\theta}\in K_0}\omega(\tilde{\theta}),0\big) + \nu\big\}$, which is compact and satisfies $K_v\subset K_1$ due to the fact that $\nu<1$. We can now use the exact same argument of Step 4 in the proof of Theorem \ref{thm:closeness_of_trajs} to establish that  \eqref{thetaboundkl} also holds for every solution starting in $K_0$ of the unrestricted HDS \eqref{eq:ms_orig_hybrid_sys}. \hfill $\blacksquare$ 

\vspace{-0.2cm}
\subsection{Proofs of Propositions \ref{prop:switching_topology_es}-\ref{prop:global_es_example}} 
\label{proofsPropositions}
\noindent 
\textbf{Proof of Proposition \ref{prop:switching_topology_es}:} Since $\mathcal{Q}_{u}=\emptyset$, $\eta_2>0$ and $T_{\circ}\geq 0$ can be arbitrary. \thref{asmp:ex_synchronization_A1} guarantees that the subset $\mathcal{S}$ is UAS for the averaged vector field (\ref{asmp:ex_synchronization_A1}) for each fixed $z_1\in\mathcal{Q}$. 
Hence, since system \eqref{flowsautomatonswitching} satisfies Assumption \ref{asmp:regularity}, and $\mathcal{A}$ is compact, by \cite[Thm. 7.12 and Lemma 7.20]{bookHDS}, the set $\mathcal{A}$ is SGPAS as $\eta_1\rightarrow 0^+$ with respect to some basin of attraction $\mathcal{B}_{\mathcal{A}}$. Therefore, by \thref{thm:avg_UAS_implies_org_PUAS}, we conclude that $\mathcal{A}$ is SGPAS as $(\varepsilon,\eta_1)\rightarrow 0^+$ for \eqref{eq:ms_orig_hybrid_sys} with respect to $\mathcal{B}_{\mathcal{A}}$. \hfill $\blacksquare$

\vspace{0.2cm}
\noindent 
\textbf{Proof of Proposition \ref{prop:ex2_es_intermittent}:} We show that the average HDS \eqref{eq:switched_averaged_sys} renders UGAS the set $\mathcal{A}$, such that Theorem \ref{thm:avg_UAS_implies_org_PUAS} can be directly used. Indeed, consider the Lyapunov function $V(\bar{x}) = J(\bar{x})-J(x^*)$, and note that for each $z_1\in\mathcal{Q}$, we have $\nabla V(\bar{x})^\intercal \bar{f}_{z_1}(\bar{x}) = (2-z_1)\nabla J(\bar{x})^\intercal P(\bar{x})\nabla J(\bar{x})$. Using strong convexity and globally Lipschitz of $\nabla J$, we obtain: For~$z_1=1$: $\nabla V(\bar{x})^\intercal \bar{f}_1(\bar{x})\leq \frac{2M_P L_J}{\mu}V(\bar{x})$; For~$z_1=2$: $\nabla V(\bar{x})^\intercal \bar{f}_2(\bar{x})\leq \frac{2M_P L_J}{\mu}V(\bar{x})$; For~$z_1=3$:  $\nabla V(\bar{x})^\intercal \bar{f}_3(\bar{x})\leq -2\lambda_P \mu V(\bar{x})$. By \cite[Proposition 3]{PoTe17Auto}, we obtain that for $\eta_1>0$ and for $\eta^*_2>0$ sufficiently small, the set $\mathcal{A}=\{x^{\star}\}\times C_z$ is UGAS for the average HDS \eqref{eq:switched_averaged_sys}. Hence, by Theorem \ref{thm:avg_UAS_implies_org_PUAS}, $\mathcal{A}$ is SGPAS as $\varepsilon\to 0^+$ for the original system \eqref{eq:ms_orig_hybrid_sys}. \hfill  $\blacksquare$

\vspace{0.2cm}
\noindent 
\textbf{Proof of Proposition \ref{prop:global_es_example}:} We show that $\mathcal{H}_2^{\text{ave}}$ renders UGAS the set $\mathcal{A}$. To do this, we consider the Lyapunov function $V(z)=\tilde{J}_{z}(x)-J(x^{\star})$, which during flows satisfies $\dot{V}= -\sum_{i=1}^r\langle b_i(\bar{x}),\nabla \tilde{J}_{\bar{z}}(\bar{x})\rangle^2 \leq 0,~~\forall~(\bar{x},\bar{z})\in C$. If $\dot{V}=0$ for some $(\bar{x},\bar{z})\in C$, then we have three possible cases:
\begin{enumerate}
    \item $\nabla \tilde{J}_{\bar{z}}(\bar{x})\neq 0$, but $P(\bar{x})[\nabla \tilde{J}_{\bar{z}}(\bar{x})]=0$, which implies that $\nabla \tilde{J}_{\bar{z}}(\bar{x})\in (T_{\bar{x}} M)^\perp$ due to item (a) in \thref{asmp:mfd_ex_A2}. Consequently, we must have $\nabla_M \tilde{J}_{\bar{z}}(\bar{x}) = 0$; or 
    \item $\nabla \tilde{J}_{\bar{z}}(\bar{x}) = 0$, which again implies that $\nabla_M \tilde{J}_z(\bar{x}) = 0$; or
    \item $\nabla \tilde{J}_{\bar{z}}(\bar{x}) \neq 0$ and $P(\bar{x})[\nabla \tilde{J}_{\bar{z}}(\bar{x})]\neq 0$, but $\langle \nabla \tilde{J}_z(\bar{x}),P(\bar{x})[\nabla \tilde{J}_{\bar{z}}(\bar{x})] \rangle = 0.$
    In that case, we can decompose $\nabla \tilde{J}_{\bar{z}}(\bar{x}) = \nabla_M \tilde{J}_{\bar{z}}(\bar{x}) + \nabla_M^{\perp}\tilde{J}_{\bar{z}}(\bar{x})$ where $\nabla_M^{\perp}\tilde{J}_{\bar{z}}(\bar{x})\in (T_x M)^\perp$, and therefore $P(\bar{x})[\nabla \tilde{J}_{\bar{z}}(\bar{x})]$ $= P(\bar{x})[\nabla_M \tilde{J}_{\bar{z}}(\bar{x})]$, and 
    $\langle \nabla_M \tilde{J}_z(\bar{x}),P(\bar{x})[\nabla_M \tilde{J}_{\bar{z}}(\bar{x})]\rangle = 0.$
    However, by definition of $P$, we have
    $\sum_{i=1}^r\langle b_i(\bar{x}),\nabla_M \tilde{J}_{\bar{z}}(\bar{x})\rangle^2 = 0 \implies P(\bar{x})[\nabla_M \tilde{J}_{\bar{z}}(\bar{x})] = 0,$
    which contradicts the assumption that $P(\bar{x})[\nabla \tilde{J}_{\bar{z}}(\bar{x})] = P(\bar{x})[\nabla_M \tilde{J}_{\bar{z}}(\bar{x})]\neq 0$. Hence, this case cannot happen.
\end{enumerate}
Therefore, $\dot{V}=0$ for some $(\bar{x},\bar{z})\in C$ implies that $\nabla_M \tilde{J}_{\bar{z}}(\bar{x}) = 0$. Due to item (b) in \thref{asmp:mfd_ex_A2} and item (b) in \thref{asmp:mfd_ex_A1}, and since flows are allowed only in the set $C$, which contains no critical points other than in $\mathcal{A}$, it follows that $\dot{V}=0$ if and only if $(\bar{x},\bar{z})\in\mathcal{A}$. 
Next, observe that immediately after jumps we have for all $(\bar{x},\bar{z})\in D$: $\Delta V = V(\bar{x}^+,\bar{z}^{+})-V(\bar{x},\bar{z}) = \tilde{J}_{\bar{z}^+}(\bar{x}) - \tilde{J}_{\bar{z}}(\bar{x}) \leq -\delta < 0$, since, by definition of the jump map we have $\bar{z}^+\in\{z\in\mathcal{Q}: \tilde{J}_{z}(\bar{x})=\min_{\tilde{z}\in\mathcal{Q}}\tilde{J}_{\tilde{z}}(\bar{x})\}$ and, using the structure of $D$, 
\begin{align*}
    (\bar{x},\bar{z})\in D \implies \tilde{J}_{\bar{z}}(\bar{x})-\min_{\tilde{z}\in\mathcal{Q}}\tilde{J}_{\tilde{z}}(\bar{x})\geq \delta > 0. 
\end{align*}
By \cite[Theorem 3.19]{HybridControlBooks}, $\mathcal{A}$ is UGAS for the HDS \eqref{eq:ms_avg_hybrid_sys}. By \thref{corollary1}, $\mathcal{A}$ is GPAS as $\varepsilon\to0^+$ for the HDS  \eqref{eq:ms_orig_hybrid_sys}. \hfill $\blacksquare$

\section{Conclusions}
\label{seconclusions}
A second-order averaging result is introduced for a class of hybrid dynamical systems that combine differential and difference inclusions. Under regularity conditions, this result establishes the closeness of solutions between the original and average dynamics, enabling semi-global practical asymptotic stability when the average hybrid system has an asymptotically stable compact set. The findings are demonstrated through the analysis and design of various hybrid and highly oscillatory algorithms for model-free control and optimization problems.  {The theoretical tools developed in this paper also enable the study of systems for which global control Lyapunov functions, as defined in \cite{scheinker2017model}, do not exist. For additional applications, we refer the reader to the recent manuscript \cite{MahmoudPoveda2024CLF}. 
\tcb{Other future research directions will study incorporating resets into the controllers \cite{li2012extremum,PovedaNaliAuto20},  stochastic phenomena, the development of hybrid source-seeking controllers with multi-obstacle avoidance capabilities \cite{PatentBenosmanPoveda}, as well as developing experimental validations of the proposed algorithms.}}

\bibliographystyle{plain}

\bibliography{References,Ref1}

\appendix

\section{Proofs of Auxiliary Lemmas \ref{lemma0technical}-\ref{keytechnicallemma1}}
\textbf{Proof of Lemma \ref{lemma0technical}:} By \thref{asmp:A1}, $\phi_1$ is $\mathcal{C}^1$ in $x$ and $\tau_1$ and $\mathcal{C}^0$ in $z$, and therefore so is $u_1$ since it is the integral of $\phi_1$ with respect to $\tau_2$. Moreover, $\phi_2$ is $\mathcal{C}^0$ in $x$, $\tau_1$, and $z$. It follows that $\psi_m$ is $\mathcal{C}^0$ in $x$, $\tau_1$, and $z$, since it is the integral with respect to $\tau_2$ of terms that involve $\partial_x \phi_1,\,\partial_x u_1,\,\phi_2,\, \phi_1,$ and $u_1$, which are all $\mathcal{C}^0$ in $x$, $\tau_1$, and $z$. Consequently, $\bar{h}$ is $\mathcal{C}^0$ in all its arguments. In addition, since $\bar{f}$ is the integral with respect to $\tau_1$ of $\bar{h}$, it follows that $\bar{f}$ is also $\mathcal{C}^0$ in all its arguments. The conclusion of the lemma follows from the fact that the set $K$ is compact and that all functions are periodic in $\tau_1$, and therefore $\bar{h}$ is uniformly Lipschitz continuous and bounded in $(K\cap C_{\delta})\times\mathbb{R}_{\geq 0}$, and the vector field $\bar{f}$ is uniformly Lipschitz continuous on $K\cap C_{\delta}$. \hfill $\blacksquare$ 

\textbf{Proof of Lemma \ref{keylemma1periodic}:} By \thref{asmp:A1}, $\phi_1$ is continuous with respect to $\tau_2$, and is $\mathcal{C}^1$ with respect to $x$ and $\tau_1$, i.e. continuously differentiable and its differential (the jacobian matrix $[\partial_x \phi_1\,\,\partial_{\tau_1}\phi_1]$) is locally Lipschitz continuous. Since $u_1$ is the integral of $\phi_1$ with respect to $\tau_2$, it follows that $u_1$, $\partial_x u_1$, and $\partial_{\tau_1}u_1$ are $\mathcal{C}^0$ in all arguments. In addition, we note that $\phi_2$ is $\mathcal{C}^0$ in $x$ and $\tau_1$, and continuous with respect to $\tau_2$. Therefore, $u_2$ is the integral with respect to $\tau_2$ of terms that are locally Lipschitz continuous in $x$ and $\tau_1$, and continuous in $\tau_2$. Consequently, $u_2$ is also $\mathcal{C}^0$ in all arguments. Moreover, since all terms are periodic in $\tau_1$ and $\tau_2$, we obtain that all terms are globally Lipschitz in $\tau_1$ and $\tau_2$. This establishes the inequalities for the maps $\sigma_k$ for $k\in\{1,2\}$. An identical argument establishes the inequality for the map $\bar{\sigma}.$ \hfill $\blacksquare$

\textbf{Proof of Lemma \ref{keytechnicallemma1}:} By invoking \cite[Lemma 2]{TeelNesicAveraging}, there exist functions $\tilde{\sigma}_k,\tilde{\sigma}$, constructed using the saturation function as in \cite[Lemma 2]{TeelNesicAveraging}, such that, due to Lemma \ref{keylemma1periodic}, the properties in items (a), (b), and (c) hold. To establish item (d), we note that since $\tilde{\sigma}_1=\sigma_1$ in the set $C|_{K,\delta}$, and by using the construction of $\sigma_1$ in \eqref{constructionsigma1} and the definition of $u_1$ in \eqref{eq:u1def}, we directly obtain $\partial_{\tau_2}u_1=\phi_1$, which establishes \eqref{importantPropLemma2}. Similarly, since
\begin{subequations}\label{chainsformulae}
\begin{align}
\partial_{\tau_2}\widetilde{{\sigma}}_2&=\phi_2+\psi_m+\psi_p-\bar{h}-\partial_{\tau_1}u_1-\partial_{\tau_2}(\partial_xu_1\cdot u_1),\\
\partial_{x}\widetilde{{\sigma}}_1\,\phi_1&=\partial_x u_1\cdot \phi_1,~~~~\partial_{\tau_1}\widetilde{{\sigma}}_1=\partial_{\tau_1}u_1,
\end{align}
\end{subequations}
and since $\partial_{\tau_2}(\partial_xu_1\cdot u_1)=\partial_{\tau_2} (\partial_x u_1)\cdot u_1+ \partial_xu_1\cdot \phi_1$, and

\vspace{-0.2cm}
\begin{small}
\begin{align}
\partial_{\tau_2}\left(\partial_{x}\int_0^{\tau_2} \phi_1(\bar{\theta},\tau_1,s_2)\,ds_2\right)&=\partial_{\tau_2}\left(\int_0^{\tau_2} \partial_{x}\phi_1(\bar{\theta},\tau_1,s_2)\,ds_2\right)\notag\\
&=\partial_x \phi_1,\label{lastchain1}
\end{align}
\end{small}
we can combine \eqref{chainsformulae}-\eqref{lastchain1} to obtain \eqref{important2PropLemma2}. To establish items (e) and (f), let $M_k>0$ satisfy $|\phi_k(\theta,\tau)|\leq M_k$, for $k\in\{1,2\}$, for all $(\theta,\tau)\in C|_{K,\delta}\times\mathbb{R}^{2}_{\geq0}$, which exists due to the continuity of $\phi_k$, the compactness of $C|_{K,\delta}$ and the periodicity of $\phi_k$ with respect to $\tau$. In addition, since $\Phi$ is LB, there exist $M_{z}>0$ such that $\Phi(C|_{K,\delta})\subset M_{z}\mathbb{B}$. Moreover, since  $\widetilde{{\sigma}}_1$ and $\widetilde{{\sigma}}_2$ are globally Lipschitz by construction, it follows that the generalized Jacobians $\partial_z\widetilde{\sigma}_1$ and $\partial_x\widetilde{\sigma}_2$ are OSC, compact, convex, and bounded \cite{clarke2008nonsmooth}. Therefore, $\forall (\theta,\tau)\in C|_{K,\delta}\times\mathbb{R}_{\geq0}^2$, $\forall a_3\in\alpha_3(\theta,\tau)$, and $\forall a_4\in\alpha_4(\theta,\tau)$, we have 
\begin{align*}
    \lvert a_3\lvert &\leq L_2 M_1+ L_1 M_{z} + L_1 M_2  + L_1=:M_3, \\ \lvert a_4\lvert &\leq L_2 M_2 + L_2 M_{z}=:M_4,
\end{align*}
which establishes the result. 
\hfill $\blacksquare$

\section{Construction of an adversarial input on $\mathbb{S}^n$}
Let $N\in\mathbb{N}_{\geq 1}$, and let $\{x_1,\dots,x_N\}\in\mathbb{S}^n$ be a set of points such that there exists $\Delta\in\mathbb{R}_{>0}$ such that, if $\mathcal{B}_i=\mathbb{S}^n\cap(\{x_i\}+\Delta\mathbb{B})$, then $\mathcal{B}_i\cap\mathcal{B}_j=\emptyset$ for all $i\neq j$. Consider the auxiliary functions $\chi_j:\mathbb{R}\rightarrow \mathbb{R}_{\geq 0}$, given by
\begin{align*}
    \chi_1(r)&:= \begin{cases}
        \exp\left(- r^{-1}\right) & r>0\\
        0 & r\leq 0
    \end{cases}, \\
    \chi_2(r)&:= \frac{\chi_1(r)}{\chi_1(r) + \chi_1(1-r)}.
\end{align*}
Using the function $\chi_2$, the constant $\Delta$, and an arbitrary choice of $\epsilon\in(0,1)$, we define the smooth bump functions $\varphi_i:\mathbb{S}^n\rightarrow [0,1]$, for $i\in\{1,\dots,N\}$, by the expression:
\begin{align}\label{bumpfunction}
    \varphi_i(x)&= 
        1-\chi_2\left(\frac{1-\langle x,x_i\rangle-(1-\epsilon)\Delta}{\epsilon\Delta}\right).
\end{align}
It is straightforward to verify that $\varphi_i(x)=0$, for all $x\not\in\mathcal{B}_i$. In other words, $\text{supp}(\varphi_i)\subset \mathcal{B}_i$, where $\text{supp}(\varphi_i)=\{x\in\mathbb{S}^n:\,\varphi_i(x)\neq 0\}$. Moreover, it can be shown that there exists $\mathcal{U}_i\subset \text{supp}(\varphi_i)$ such that $\varphi_i(x)=1$ for all $x\in\mathcal{U}_i$. Using the above construction, let $f_e:\mathbb{S}^n\ni x\mapsto T_x\mathbb{S}^n$ be given by
\begin{align*}
    f_e(x):=\delta_1\sum_{i=1}^N {\varphi_i(x)(x_i-\langle x_i,x\rangle x)},
\end{align*}
where $\delta_1>0$ is a tuning parameter. It is straightforward to see that the map $f_e$ is smooth and that 
\begin{align*}
    \sup_{x\in\mathbb{S}^n}|f_e(x)|^2\leq\delta_1^2\max_{i\in\{1,\dots,N\}}\sup_{\tilde{x}\in\mathcal{B}_i}(1-\langle x_i,\tilde{x}\rangle^2).
\end{align*}
Define the candidate Lyapunov functions $V_i:\mathbb{S}^n\rightarrow\mathbb{R}_{\geq 0}$ by $V_i(x):=1-\langle x_i,x\rangle$, then observe that, due to the properties of the functions $\varphi_i$, the derivatives of the functions $V_i$ along the vector field $f_e$ are given by
\begin{align*}
    \left\langle\nabla_{\mathbb{S}^n}V_i(x),f_e(x)\right\rangle < -{\delta_1(1-\langle x_i,x\rangle^2)}\leq 0,
\end{align*}
for all $x\in\mathcal{U}_i$, and that $\left\langle\nabla_{\mathbb{S}^n}V_i(x),f_e(x)\right\rangle\iff x = x_i$. In other words, the vector field $\dot{x}= f_e(x)$, locally stabilizes all the points in the set $\{x_1,\dots,x_N\}\in\mathbb{S}^n$. Moreover, if $f_i:\mathbb{S}^n \ni x\rightarrow T_x\mathbb{S}^n$ is another smooth vector field such that 
\begin{align*}
\left\langle\nabla_{\mathbb{S}^n}V_i(x),f_i(x)\right\rangle<\delta_3(1-\langle x_i,x\rangle^2),
\end{align*} 
for some $\delta_3>0$ and for all $x\in\mathcal{U}_i$, then we have that:
\begin{align*}
\left\langle\nabla_{\mathbb{S}^n}V_i(x),f_e(x)+f_i(x)\right\rangle < -(\delta_1-\delta_3)(1-\langle x_i,x\rangle^2).
\end{align*}
That is, if $\delta_3<\delta_1$, then the perturbation $f_i(x)$ does not destroy the (local) stability of the the point $x_i$. 
\end{document}